# ABSOLUTE CONTINUITY OF SYMMETRIC MARKOV PROCESSES

By Z.-Q. Chen[1], P. J. Fitzsimmons, M. Takeda, J. Ying and T.-S. Zhang[2]

*University of Washington, University of California, San Diego, Tohoku University, Fudan University and University of Manchester*

We study Girsanov's theorem in the context of symmetric Markov processes, extending earlier work of Fukushima–Takeda and Fitzsimmons on Girsanov transformations of "gradient type." We investigate the most general Girsanov transformation leading to another symmetric Markov process. This investigation requires an extension of the forward–backward martingale method of Lyons–Zheng, to cover the case of processes with jumps.

**1. Introduction and preliminaries.** Our aim in this work is to study change-of-measure phenomena ("Girsanov" theorems) for general symmetric Markov processes. Our results extend both the earlier work of Fukushima and Takeda [9] and that of Fitzsimmons [5] (who was concerned only with symmetric diffusions). Our work also contains Theorem 2.7 of [12] as a special case.

Before setting down the precise context in which we shall be working, let us briefly describe our results. Let $X = (X_t)$ be a symmetric (i.e., reversible) Markov process, with symmetry measure $m$, state space $E$ and distribution $\mathbb{P}_x$ when started in state $x \in E$. Given a strictly positive element $\rho$ of the Dirichlet space of $X$, let $M^\rho$ be the martingale part in the Fukushima decomposition of $\rho(X_t) - \rho(X_0)$, define a local martingale $M$ by the formula $M_t := \int_0^t [\rho(X_{s-})]^{-1} dM_s^\rho$ and let $Z^\rho$ be the "stochastic exponential" of $M$; that is, $Z^\rho$ is the unique solution of $Z_t^\rho = 1 + \int_0^t Z_{s-}^\rho dM_s$. Then $Z^\rho$ is a

Received June 2002; revised July 2003.

[1]Supported in part by NSF Grant DMS-00-71486.

[2]Supported by the Norwegian Research Council.

positive supermartingale, and so determines a new family $(\widehat{\mathbb{P}}_x)_{x \in E}$ of probability measures governing a second symmetric Markov process $\widehat{X}$ on $E$, with symmetry measure $\mu(dx) := \rho(x)^2 \, m(dx)$.

Suppose, for example, that $E = \mathbb{R}^d$ and that the (nonpositive definite) infinitesimal generator $L$ of $X$ has the form

$$(1.1) \quad Lf(x) = L^c f(x) + \int_{\mathbb{R}^d} [f(y) - f(x)] N(x, dy) - k(x) f(x), \qquad x \in \mathbb{R}^d,$$

where $L^c$ is a second-order differential operator with $L^c 1 = 0$, $N$ is a kernel on $\mathbb{R}^d$ and $k \geq 0$. Suppose also that $D(L)$, the domain of $L$, is closed under the formation of products, and that $\rho \in D(L)$ is strictly positive. Then the generator $\widehat{L}$ of the transformed process $\widehat{X}$ is given by

$$(1.2) \quad \begin{aligned} \widehat{L}f(x) &= \rho(x)^{-1}[L(\rho f)(x) - f(x) L\rho(x)] \\ &= L^c f(x) + B_\rho f(x) + \int_{\mathbb{R}^d} [f(y) - f(x)] \widehat{N}(x, dy), \end{aligned}$$

where $B_\rho f := \rho^{-1}[L^c(\rho f) - f L^c \rho - \rho L^c f]$ is a first-order differential operator and $\widehat{N}(x, dy) = [\rho(y)/\rho(x)] N(x, dy)$. It is an important but very challenging problem to characterize the domain $D(\widehat{L})$ of $\widehat{L}$. An equivalent, but more tractable, way to proceed is to characterize the bilinear *Dirichlet form* associated with $\widehat{L}$. Let $\mathcal{E}(f, g) := -\int Lf \cdot g \, dm$ and $\widehat{\mathcal{E}}(f, g) := -\int \widehat{L}f \cdot g \, d\mu$ denote the Dirichlet forms corresponding to $X$ and $\widehat{X}$. When expressed in terms of these bilinear forms, (1.1) and (1.2) become

$$(1.3) \quad \begin{aligned} \mathcal{E}(f, g) &= \mathcal{E}^c(f, g) + \int_{\mathbb{R}^d} \int_{\mathbb{R}^d} [f(y) - f(x)] \cdot [g(y) - g(x)] J(dx, dy) \\ &\quad + \int_{\mathbb{R}^d} f(x) g(x) k(x) m(dx) \end{aligned}$$

and

$$(1.4) \quad \widehat{\mathcal{E}}(f, g) = \mathcal{E}_\rho^c(f, g) + \int_{\mathbb{R}^d} \int_{\mathbb{R}^d} [f(y) - f(x)] \cdot [g(y) - g(x)] \widehat{J}(dx, dy),$$

where $\mathcal{E}^c(f, g) := -\int L^c f \cdot g \, dm$, $\mathcal{E}_\rho^c(f, g) := -\int [L^c + B_\rho] f \cdot g \, d\mu$, $J(dx, dy) := \frac{1}{2} m(dx) N(x, dy)$ and $\widehat{J}(dx, dy) := \frac{1}{2} \mu(dx) \widehat{N}(x, dy) = \rho(x) \rho(y) \times J(dx, dy)$. In the general context in which we shall be working, Dirichlet forms are more convenient objects than their associated infinitesimal generators, and formulas like (1.4) will be the main focus of our study.

The following is concrete instance in which (1.1) holds.

EXAMPLE 1.1.   Let $E = \mathbb{R}^d$, $m(dx) = dx$ (Lebesgue measure), and let $X$ be the Lévy process on $\mathbb{R}^d$ that is the sum of Brownian motion on $\mathbb{R}^d$ and an independent rotationally symmetric $\alpha$-stable process on $\mathbb{R}^d$, for some $0 < \alpha < 2$ and $d \geq 3$. Clearly $X$ is an $m$-symmetric Hunt process. Its Dirichlet form $(\mathcal{E}, \mathcal{F})$ is given by

$$\mathcal{F} = \{f \in L^2(\mathbb{R}^d, dx) : \nabla f \in L^2(\mathbb{R}^d, dx)\}$$

and

$$
\begin{aligned}
(1.5) \qquad \mathcal{E}(f, g) &= \frac{1}{2} \int_{\mathbb{R}^d} \nabla f(x) \cdot \nabla g(x) \, dx \\
&\quad + \int_{\mathbb{R}^d} \int_{\mathbb{R}^d} [f(y) - f(x)] \cdot [g(y) - g(x)] \frac{c(d, \alpha)}{|x - y|^{d + \alpha}} \, dx \, dy,
\end{aligned}
$$

where $c(d, \alpha) > 0$ is a constant depending only on $d$ and $\alpha$. In this case we have, in (1.1),

$$(1.6) \qquad L^c = \tfrac{1}{2} \Delta, \qquad N(x, dy) = \frac{c(d, \alpha)}{|x - y|^{d + \alpha}} \, dy, \qquad k(x) \equiv 0.$$

In the case at hand, (1.4) becomes

$$
\begin{aligned}
(1.7) \qquad \widehat{\mathcal{E}}(f, g) &= \frac{1}{2} \int_{\mathbb{R}^d} \nabla f(x) \cdot \nabla g(x) \, \rho(x)^2 \, dx \\
&\quad + \int_{\mathbb{R}^d} \int_{\mathbb{R}^d} [f(y) - f(x)] \cdot [g(y) - g(x)] \frac{c(d, \alpha) \rho(x) \rho(y)}{|x - y|^{d + \alpha}} \, dx \, dy.
\end{aligned}
$$

The type of change of measure considered above was studied in [9], for $\rho$ the $\alpha$-potential of a bounded strictly positive function on $E$. In Section 2, we extend and complete the work of Fukushima and Takeda by obtaining a complete description of the Dirichlet form associated with $\widehat{X}$ for an arbitrary strictly positive $\rho$ in the Dirichlet space of $X$. Our method is a modification of that found in the work of Chen and Zhang [3]. Of critical importance is Lemma 2.2 which extends the Lyons–Zheng forward–backward martingale decomposition to the context of symmetric Markov processes with jumps.

When $X$ is a diffusion, the change of measure determined by $Z^\rho$ is (modulo localization) the most general Girsanov transformation leading to another symmetric diffusion. This assertion is the principal result of [5]. The situation is more complex when $X$ has jumps. In Sections 3 and 4, we investigate the most general change of measure leading to a second symmetric Markov process $\widehat{X}$, and we take the first steps in describing the associated Dirichlet form. A formula like (1.4) holds even in this more general situation, although a zero-order term might be present, and the density $\rho(x)\rho(y)$ linking the measures $\widehat{J}$ and $J$ must be replaced by a more general symmetric function of $(x, y)$. Our results here are somewhat less comprehensive than those of Section 2 or of [5], the main unresolved difficulty being the description of a core for the Dirichlet space of $\widehat{X}$. Nevertheless, we find an explicit expression for the Dirichlet form of $\widehat{X}$ that is valid for a large class of functions in the Dirichlet space. It may be helpful for the reader to keep in mind the concrete Example 1.1 when considering the general results of this paper.

In the remainder of this section, we establish our setting and notation. Let $E$ be a topological space that is homeomorphic to a co-analytic subset of a compact metric space (a "metric co-Souslin space"), with Borel $\sigma$-algebra $\mathcal{B}(E)$. Let $m$ be a $\sigma$-finite measure on $\mathcal{B}(E)$ with supp$[m] = E$. We denote by $\mathcal{B}(E \times E)$ the product $\sigma$-algebra on $E \times E$. Let $X = (\Omega, \mathcal{M}, \mathcal{M}_t, \theta_t, X_t, \mathbb{P}_x)$ be a Borel right Markov process with state space $E$, lifetime $\zeta$, transition semigroup $(P_t)_{t>0}$ and resolvent $(U^q)_{q>0}$. In more detail, the right-continuous process $[0, +\infty[ \ni t \mapsto X_t$ is defined on the sample space $(\Omega, \mathcal{M})$, with its minimal (augmented) admissible filtration $\{\mathcal{M}_t\}_{t \geq 0}$, and under the law $\mathbb{P}_x$ is a strong Markov process with initial condition $X_0 = x$. The shift operators $\theta_t$, $t \geq 0$, satisfy $X_s \circ \theta_t = X_{s+t}$ identically for $s, t \geq 0$. Adjoined to the state space $E$ is an isolated point $\Delta \notin E$; the process $X$ retires to $\Delta$ at its "lifetime" $\zeta := \inf\{t : X_t = \Delta\}$. Throughout this paper we assume that $X$ is $m$-symmetric. More precisely, $(P_t)$ may be extended into a symmetric operator semigroup on $L^2(m)$; that is,

$$(f, P_t g) = (P_t f, g), \qquad f, g \in L^2(m),$$

where $(u, v) := \int_E u(x) v(x) \, m(dx)$ is the natural inner product in $L^2(m) := L^2(E; m)$. By the theory of Dirichlet forms, there exists a symmetric *Dirichlet form* $(\mathcal{E}, \mathcal{F})$ associated with $X$:

$$\mathcal{F} = \left\{ u \in L^2(m) : \sup_{t > 0} \frac{1}{t}(u - P_t u, u) < \infty \right\},$$

$$\mathcal{E}(u, v) = \lim_{t \downarrow 0} \frac{1}{t}(u - P_t u, v), \qquad u, v \in \mathcal{F}.$$

For any $q > 0$, set

$$\mathcal{E}_q(u, v) := \mathcal{E}(u, v) + q(u, v), \qquad u, v \in \mathcal{F}.$$

Then $\mathcal{F}$ becomes a Hilbert space with inner product $\mathcal{E}_q$ for any $q > 0$. We call the corresponding norm the $\mathcal{E}_q$-norm. In view of the work in [6] and [19], the Dirichlet form $(\mathcal{E}, \mathcal{F})$ is *quasi-regular*. Thus, by Theorem 3.7 of [2], the process is quasi-homeomorphic to the Hunt process associated with a regular Dirichlet form on a locally compact separable metric space, so all of the results of [8] (established there for regular Dirichlet forms) apply to $X$ and its Dirichlet form. For the reader uninterested in applications to infinite-dimensional settings, it is safe to assume from now on that $X$ is the Hunt process on a locally compact separable metric state space (such as a Euclidean domain) associated with a regular Dirichlet form. See Chapter VI of [19] for more on the "transfer method," by which the quasi-regular case is reduced to the regular case.

Let $A = (A_t)$ be any increasing additive functional (AF), by which we mean that besides additivity and right-continuity we also assume that $0 \le A_t < \infty$ for $t < \zeta$. We can (and do) take its perfected version. We use $f * A$ to denote the functional

$$(f * A)_t := \int_0^t f(X_s) \, dA_s,$$

if $f$ is a Borel function on $E$, and

$$(F * A)_t := \int_0^t F(X_{s-}, X_s) \, dA_s,$$

if $F$ is a Borel function on $E \times E$. The *bivariate Revuz measure* $\nu_A$ of $A$ (computed with respect to $X$ and $m$) is defined for $F \in p\mathcal{B}(E \times E)$ by

$$\nu_A(F) := \uparrow \lim_{t \downarrow 0} \frac{1}{t} \mathbb{E}_m(F * A)_t.$$

The usual Revuz measure $\mu_A$ of $A$ is nothing but the second marginal measure of $\nu_A$, namely, $\mu_A(dx) = \nu_A(1 \otimes dx)$. The mapping $A \mapsto \mu_A$ establishes a one-to-one correspondence between the class of positive continuous additive functionals (PCAFs) of $X$ and the class of smooth measures of $(\mathcal{E}, \mathcal{F})$, and is usually known as the Revuz correspondence.

A well-known consequence of symmetry is that, for q.e. $x \in E$,

$$\mathbb{P}_x(\{\omega \in \Omega : X_{t-}(\omega) \text{ exists in } E \text{ for all } t < \zeta\}) = 1.$$

Without any real loss of generality, we assume the exceptional set [of those $x \in E$ for which (1.5) fails] to be empty. Adjoin the cemetery $\Delta$ to $E$ as an isolated point to form $E_\Delta$; the Borel $\sigma$-field on $E_\Delta$ is denoted $\mathcal{B}(E_\Delta)$. The jump behavior of $X$ is described by a pair $(N, H)$, *the Lévy system* of $X$, in which $N$ is a kernel from $(E, \mathcal{B}(E))$ to $(E_\Delta, \mathcal{B}(E_\Delta))$ satisfying $N(x, \{x\}) = 0$ for any $x \in E$, and $H$ is a PCAF of $X$ with bounded one-potential, such that for any measurable function $F \in p\mathcal{B}(E \times E)$, the dual predictable projection (or compensator) of the homogeneous random measure

$$\eta(\omega, dt) := \sum_{s>0} F(X_{s-}(\omega), X_s(\omega)) \mathbb{1}_{\{X_{s-}(\omega) \ne X_s(\omega)\}} \varepsilon_s(dt)$$

is $NF * H$, where $NF(x) := \int_{E_\Delta} N(x, dy) F(x, y)$. (Here $\varepsilon_s$ is the unit point mass at $s$.) The special case (1.1) occurs when $H_t \equiv t$. Set $J(dx, dy) := \frac{1}{2} \mathbb{1}_E(y) \mu_H(dx) N(x, dy)$, the jump measure for $X$. Note that $J$ is invariant with respect to the mapping $(x, y) \to (y, x)$.

We will use $\zeta_p$ and $\zeta_i$ to denote, respectively, the predictable and totally inaccessible parts of the lifetime $\zeta$. Let $\kappa(dx) := N(x, \{\Delta\}) \mu_H(dx)$ be the Revuz measure of the compensator $A^\kappa$ of the AF $\mathbb{1}_{[\![\zeta_i, \infty[\![}$; $\kappa$ is the killing measure for $X$.

Each $u \in \mathcal{F}$ admits a *quasi-continuous* $m$-version $\tilde{u}$ for which the process $t \mapsto \tilde{u}(X_t)$ is right-continuous on $[0, \infty[$ with left limits [equal to $\tilde{u}(X_{t-})$ for $t \neq \zeta_p$] on $]0, \infty[$, $\mathbb{P}_x$-a.s. for q.e. $x \in E$. For such $u \in \mathcal{F}$, we have Fukushima's decomposition (**(year?)**, Theorem 5.2.2)

$$\tilde{u}(X_t) - \tilde{u}(X_0) = M_t^u + N_t^u,$$

where $M^u$ is a martingale additive functional (MAF) of finite energy and $N^u$ is a continuous additive functional (CAF) of zero energy. Moreover, $M^u$ may be further decomposed into the sum of a continuous part and a purely discontinuous part

$$M^u = M^{u,c} + M^{u,d}.$$

Let $\mathcal{E}^c(u, u) := \frac{1}{2}\mu_{\langle M^{u,c} \rangle}(E)$. Then $\mathcal{E}$ admits a Beurling–Deny decomposition as in ([8], Theorem 5.3.1):

$$\mathcal{E}(u, u) = \mathcal{E}^c(u, u) + \int_{E \times E} (u(y) - u(x))^2 J(dx, dy) + \kappa(u^2), \qquad u \in \mathcal{F},$$

where the "diffusion" piece $\mathcal{E}^c$ is *strongly local* in the sense that $\mathcal{E}(u, v) = 0$ whenever $u, v \in \mathcal{F}$ and $u$ is constant $m$-a.e. on a neighborhood of the support of the measure $|v| \cdot m$. Here and in the sequel, we always take $u \in \mathcal{F}$ to be represented by its quasi-continuous version, and we usually drop the "tilde" from our notation.

**Notation and conventions.** The notation ":=" should be read "is defined to be." For a class $\mathcal{F}$ of functions, we use $b\mathcal{F}$ [resp. $p\mathcal{F}$ (or $\mathcal{F}^+$)] to denote the set of bounded (resp. nonnegative) functions in $\mathcal{F}$. We use both "nonnegative" and "positive" to mean $\geq 0$, and "strictly positive" to mean $> 0$. For a measure $\mu$ and a function $f$, $\mu(f) := \int f \, d\mu$. We sometimes write $L^p$ or $L^p(m)$ for $L^p(E, m)$, and $(\cdot, \cdot)$ for the inner product in $L^2(m)$. For $f, g \in \mathcal{B}(E)$, $f \otimes g(x, y) := f(x)g(y)$ and $\bar{f}(x, y) := f(y) - f(x)$ for $x, y \in E$. For a right-continuous process $H_t$ of finite variation on bounded intervals, we use $H^p$ to denote its dual predictable projection and $H^{\sim}$ to denote $H - H^p$, all computed with respect to $(X, \mathbb{P}_x, x \in E)$. The jump $M_t - M_{t-}$ will be abbreviated as $\Delta M_t$. The first *hitting time* of a set $G$ is denoted $T_G := \inf\{t > 0 : X_t \in G\}$. A hitting time is an example of a *terminal time*, which is a stopping time $T$ such that $t + T \circ \theta_t = T$ on $\{t < T\}$.

**2. Girsanov transform by multiplicative functional related to $M^{\log \rho}$.** In this section, we study Girsanov transforms of a type investigated earlier by Fukushima and Takeda [9]; our results extend and complete the work found there. Our method is a modification of that found in Chen and Zhang [3].

Throughout this section, $\rho$ is a nonnegative element of $\mathcal{F}$. We can (and do) assume that $\rho$ is quasi-continuous, and we assume that $\rho > 0$ q.e. on $E$. [Otherwise, we would deal with the part process $X$ killed upon leaving $\{x : \rho(x) > 0\}$.] We will use the convention that any function defined on $E$ is extended to be zero at the cemetery point $\Delta$; in particular, $\rho(\Delta) = 0$. By Fukushima's decomposition,

$$\rho(X_t) - \rho(X_0) = M_t^\rho + N_t^\rho, \qquad \mathbb{P}_x\text{-a.s. for q.e. } x \in E,$$

where $M^\rho$ is a square-integrable martingale AF and $N^\rho$ is a CAF of zero energy. Note that $s \mapsto \rho(X_{s \wedge \zeta_i -})$ is strictly positive and left-continuous on $]0, \zeta_p[$. Define a local martingale $M$ on the random time interval $[\![0, \zeta_p[\![$ by

$$(2.1) \qquad M_t = \int_0^{t \wedge \zeta_i} \frac{1}{\rho(X_{s-})} \, dM_s^\rho \qquad 0 \le t < \zeta_p.$$

Note that for $t < \zeta_p$,

$$\begin{aligned}
\Delta M_t &= \frac{1}{\rho(X_{t \wedge \zeta_i -})}(M_{t \wedge \zeta_i}^\rho - M_{t \wedge \zeta_i -}^\rho) \\
&= \frac{1}{\rho(X_{t \wedge \zeta_i -})}(\rho(X_{t \wedge \zeta_i}) - \rho(X_{t \wedge \zeta_i -})) \\
&= \frac{\rho(X_{t \wedge \zeta_i})}{\rho(X_{t \wedge \zeta_i -})} - 1.
\end{aligned}$$

The solution $Z_t^\rho$ of the SDE

$$(2.2) \qquad Z_t^\rho = 1 + \int_0^t Z_{s-}^\rho \, dM_s, \qquad 0 \le t < \zeta_p,$$

is a positive local martingale on the random time interval $[\![0, \zeta_p[\![$, hence a positive supermartingale. Consequently, the formula

$$d\widehat{\mathbb{P}}_x = Z_t^\rho \, d\mathbb{P}_x \qquad \text{on } \mathcal{M}_t \cap \{t < \zeta_p\} \text{ for } x \in E,$$

uniquely determines a family of probability measures on $(\Omega, \mathcal{M}_\infty)$. It is known that under these new measures, $X$ is a right Markov process on $E$; see [24], Section 62. We will use $(\widehat{X}, \mathcal{M}, \mathcal{M}_t, \widehat{\mathbb{P}}_x, x \in E)$ to denote the transformed process. Here $\widehat{X}_t(\omega) = X_t(\omega)$ but we use $\widehat{X}_t$ for emphasis when working with $\widehat{\mathbb{P}}_x$. Let $\widehat{P}_t$ be the semigroup of $\widehat{X}$, that is,

$$\widehat{P}_t f(x) = \widehat{\mathbb{E}}_x[f(\widehat{X}_t)] = \mathbb{E}_x[Z_t^\rho f(X_t)] = \mathbb{E}_x[Z_t^\rho f(X_t); t < \zeta].$$

(These transition operators need not preserve Borel measurability; this minor inconvenience can be dealt with as in Corollary 3.23 of [6].)

Before stating the next result, let us recall the definition of time-reversal operator $r_t$ on the path space. Given a path $\omega \in \{t < \zeta\}$, define

$$r_t(\omega)(s) := \begin{cases} \omega((t-s)-), & \text{for} \quad 0 \leq s < t, \\ \omega(0), & \text{for} \quad s \geq t. \end{cases}$$

Here for $r > 0$, $\omega(r-) := \lim_{s \uparrow r} \omega(s)$. It is known (see Lemma 4.1.2 of [8]) that the mapping $r_t$ preserves the measure $\mathbb{P}_m$ on $\mathcal{M}_t \cap \{t < \zeta\}$.

DEFINITION 2.1.   (i) A CAF $A_t$ is called *even* if $A_t \circ r_t = A_t$ for all $t < \zeta$.

(ii) An $m$-measurable function $u : E \to \mathbb{R}$ is *locally in $\mathcal{F}$* $(u \in \overset{\bullet}{\mathcal{F}}_{\text{loc}})$ provided there is a nest $\{G_n\}$ of finely open sets and a sequence $\{u_n\} \subset \mathcal{F}$ such that $u = u_n$, $m$-a.e. on $G_n$ for each $n$.

We recall from [4]; Theorem 2.1, that a CAF with paths locally of bounded variation (or merely of zero quadratic variation) is necessarily even.

Clearly each $u \in \overset{\bullet}{\mathcal{F}}_{\text{loc}}$ has a quasi-continuous version, and for such $u$ the continuous local martingale AF $M^{u,c}$ is well defined by

$$(2.3) \qquad M_t^{u,c} := M_t^{u_n,c} \qquad \text{for } t \leq T_{E_\Delta \setminus G_n}, n = 1, 2, \ldots.$$

The following can be regarded as an extension to functions in $\overset{\bullet}{\mathcal{F}}_{\text{loc}}$ of both Fukushima's decomposition and the Lyons–Zheng forward–backward martingale decomposition.

LEMMA 2.2.   *For $u \in \overset{\bullet}{\mathcal{F}}_{\text{loc}}$ and $t > 0$, $\mathbb{P}_m$-a.s. on $\{t < \zeta\}$,*

$$u(X_t) - u(X_0) = \tfrac{1}{2}(M_t^{u,c} - M_t^{u,c} \circ r_t)$$
$$+ \lim_{\varepsilon \downarrow 0} \sum_{0 < s \leq t} (u(X_s) - u(X_{s-})) \mathbb{1}_{\{|u(X_s) - u(X_{s-})| > \varepsilon\}}.$$

*The limit above exists in the sense of convergence in probability under $\mathbb{P}_x$, for $m$-a.e. $x \in E$.*

PROOF.   Note that when $u \in \mathcal{F}$, the martingale part $M_t^u$ in Fukushima's decomposition can be decomposed as

$$M_t^u = M_t^{u,c} + M_t^{u,j} + M_t^{u,k},$$

where $M_t^{u,c}$ is the continuous part of martingale $M^u$, and

$$M_t^{u,j} = \lim_{\varepsilon \downarrow 0} \Big\{ \sum_{0 < s \leq t} (u(X_s) - u(X_{s-})) \mathbb{1}_{\{|u(X_s) - u(X_{s-})| > \varepsilon\}} \mathbb{1}_{\{s < \zeta\}}$$
$$- \int_0^t \Big( \int_{\{y \in E \, : \, |u(y) - u(X_s)| > \varepsilon\}} (u(y) - u(X_s)) N(X_s, dy) \Big) dH_s \Big\},$$

(2.4)
$$M_t^{u,k} = \int_0^t u(X_s) N(X_s, \Delta) \, dH_s - u(X_{\zeta-}) \mathbb{1}_{\{t \geq \zeta_i\}}$$
$$= \int_0^t u(X_s) \, dA_s^\kappa - u(X_{\zeta-}) \mathbb{1}_{\{t \geq \zeta_i\}},$$

are the jump and killing parts $M^u$, respectively. See [8], Theorem A.3.9. The limit in the expression for $M^{u,j}$ is in the sense of convergence in the norm of the space of square-integrable martingales and convergence in probability under $\mathbb{P}_x$ for $m$-a.e. $x \in E$ (see [8]). So $\mathbb{P}_m$-a.s. on $\{t < \zeta\}$ it follows that

$$u(X_t) - u(X_0) = \tfrac{1}{2}(M_t^u - M_t^u \circ r_t)$$
$$= \tfrac{1}{2}(M_t^{u,c} - M_t^{u,c} \circ r_t)$$
$$+ \lim_{\varepsilon \downarrow 0} \sum_{0 < s \leq t} (u(X_s) - u(X_{s-})) \mathbb{1}_{\{|u(X_s) - u(X_{s-})| > \varepsilon\}}.$$

For $u \in \overset{\bullet}{\mathcal{F}}_{\mathrm{loc}}$, let $\{G_n\}$ be a nest of finely open sets and let $\{u_n\}$ be a sequence of functions in $\mathcal{F}$ such that $u = u_n$ q.e. on $G_n$. For each $u_n \in \mathcal{F}$, $\mathbb{P}_m$-a.s. on $\{t < \zeta\}$ we have

$$u_n(X_t) - u_n(X_0) = \tfrac{1}{2}(M_t^{u_n,c} - M_t^{u_n,c} \circ r_t)$$
$$+ \lim_{\varepsilon \downarrow 0} \sum_{0 < s \leq t} (u_n(X_s) - u_n(X_{s-})) \mathbb{1}_{\{|u_n(X_s) - u_n(X_{s-})| > \varepsilon\}}.$$

As $[0, t] \ni s \mapsto X_s$ is right-continuous with left limits in $E$ on $\{t < \zeta\}$, the lemma now follows from the above display, using (2.3) to pass to the limit as $n \to \infty$. $\square$

REMARK 2.3. By [8], Theorem A.3.9, we can also replace the indicator $\mathbb{1}_{\{|u(X_s) - u(X_{s-})| > \varepsilon\}}$ by $\mathbb{1}_{\{|e^{u(X_s)} - e^{u(X_{s-})}| > \varepsilon\}}$, and the set $\{y \in E : |u(y) - u(X_s)| > \varepsilon\}$ by $\{y \in E : |e^{u(y)} - e^{u(X_s)}| > \varepsilon\}$. Thus Lemma 2.2 can also take the following form: For $u \in \overset{\bullet}{\mathcal{F}}_{\mathrm{loc}}$,

$$u(X_t) - u(X_0) = \tfrac{1}{2}(M_t^{u,c} - M_t^{u,c} \circ r_t)$$
$$+ \lim_{\varepsilon \downarrow 0} \sum_{0 < s \leq t} (u(X_s) - u(X_{s-})) \mathbb{1}_{\{|e^{u(X_s)} - e^{u(X_{s-})}| > \varepsilon\}},$$
$$\mathbb{P}_m\text{-a.s. on } \{t < \zeta\}.$$

The convergence in the above expression is in the sense of convergence in probability under each $\mathbb{P}_x$ for $m$-a.e. $x \in E$.

LEMMA 2.4. $\widehat{P}_t$ is symmetric on $L^2(E, \rho^2 m)$.

PROOF.    Let $f, g \in b\mathcal{B}^+(E)$. By time reversal, we have

$$(\widehat{P}_t f, g)_{\rho^2 m} = \mathbb{E}_m[Z_t^\rho f(X_t) g(X_0) \rho^2(X_0)]$$
$$= \mathbb{E}_m[Z_t^\rho \circ r_t \, g(X_t) \rho^2(X_t) f(X_0)].$$

To show

$$(\widehat{P}_t f, g)_{\rho^2 m} = (f, \widehat{P}_t g)_{\rho^2 m} = \mathbb{E}_m[Z_t^\rho g(X_t) \rho^2(X_0) f(X_0)],$$

it suffices to prove the following identity:

$$(2.5) \qquad Z_t^\rho \circ r_t = Z_t^\rho \frac{\rho^2(X_0)}{\rho^2(X_t)}, \qquad \mathbb{P}_m\text{-a.s. on } \{t < \zeta\}.$$

To this end, note that by the Doléans–Dade formula ([13], Theorem 9.39), on $\{t < \zeta\}$,

$$(2.6) \begin{aligned} Z_t^\rho &= \exp\Big(M_t - \frac{1}{2}\langle M^c\rangle_t\Big) \prod_{0 < s \le t} (1 + \Delta M_s) e^{-\Delta M_s} \\ &= \exp\Big(M_t - \frac{1}{2}\langle M^c\rangle_t\Big) \prod_{0 < s \le t} \frac{\rho(X_s)}{\rho(X_{s-})} \exp\Big(1 - \frac{\rho(X_s)}{\rho(X_{s-})}\Big). \end{aligned}$$

It follows from (2.1) and (2.4) that on $\{t < \zeta\}$,

$$\begin{aligned} M_t = M_t^c + \lim_{\varepsilon \downarrow 0} \Big\{ &\sum_{0 < s \le t} \Big(\frac{\rho(X_s)}{\rho(X_{s-})} - 1\Big) \mathbb{1}_{\{|\rho(X_s) - \rho(X_{s-})| > \varepsilon\}} \\ &- \int_0^t \Big(\int_{\{y \in E_\Delta \,:\, |\rho(y) - \rho(X_s)| > \varepsilon\}} \Big(\frac{\rho(y)}{\rho(X_s)} - 1\Big) N(X_s, dy)\Big) dH_s\Big\}, \end{aligned}$$

where

$$M_t^c = \int_0^t \frac{1}{\rho(X_{s-})} \, dM_s^{\rho,c}.$$

Since $\rho > 0$ q.e. on $E$, we see that $\log \rho \in \overset{\bullet}{\mathcal{F}}_{\mathrm{loc}}$ (see [18], Corollary 6.2), and therefore, by Lemma 2.2 and Remark 2.3, we have $\mathbb{P}_m$-a.s. on $\{t < \zeta\}$,

$$(2.7) \begin{aligned} &\log \rho(X_t) - \log \rho(X_0) \\ &= \tfrac{1}{2}(M_t^c - M_t^c \circ r_t) \\ &\quad + \lim_{\varepsilon \downarrow 0} \sum_{0 < s \le t} (\log \rho(X_s) - \log \rho(X_{s-})) \mathbb{1}_{\{|\rho(X_s) - \rho(X_{s-})| > \varepsilon\}}. \end{aligned}$$

Since both $\langle M^c\rangle_t$ and $\int_0^t \int_{\{y \in E_\Delta \,:\, |\rho(y) - \rho(X_s)| > \varepsilon\}} (\rho(y)\rho(X_s)^{-1} - 1) N(X_s, dy) \, dH_s$ are even CAFs of $X$, identity (2.5) follows from (2.6) and (2.7).  □

The following result appears for symmetric diffusions as Lemma 4.4 in [5]; the proof given there is valid for general symmetric Borel right processes.

THEOREM 2.5. *If $A = (A_t)$ is a PCAF of $X$ with Revuz measure $\mu$, then the Revuz measure of $A$ as a PCAF of $\widehat{X}$ is $\rho^2 \mu$.*

In what follows, if $f \in \mathcal{B}(E)$, then we write $f \in \mathcal{L}^2(\rho \otimes \rho \cdot J)$ to mean that $\bar{f}(x, y) := f(y) - f(x)$ is square-integrable with respect to $\rho \otimes \rho \cdot J$.

THEOREM 2.6. *Let $(\widehat{\mathcal{E}}, \widehat{\mathcal{F}})$ be the symmetric Dirichlet form on $L^2(E, \rho^2 m)$ associated with $\widehat{X}$. Then*

(a)

$$\left\{ f \in \mathcal{F} : \int \rho(x)^2 \mu^c_{\langle f \rangle}(dx) < \infty \right\} \cap \mathcal{L}^2(\rho \otimes \rho \cdot J) \cap L^2(\rho^2 m) \subset \widehat{\mathcal{F}}$$

*and for $f$ in $\mathcal{F} \cap \mathcal{L}^2(\rho \otimes \rho \cdot J) \cap L^2(\rho^2 m)$ with $\int \rho(x)^2 \mu^c_{\langle f \rangle}(dx) < \infty$,*

$$
\begin{aligned}
(2.8) \qquad \widehat{\mathcal{E}}(f, f) = {} & \tfrac{1}{2} \int_E \rho(x)^2 \mu^c_{\langle f \rangle}(dx) \\
& + \int_{E \times E} (f(x) - f(y))^2 \rho(x) \rho(y) J(dx, dy);
\end{aligned}
$$

(b) *$1 \in \widehat{\mathcal{F}}$ and $\widehat{\mathcal{E}}(1, 1) = 0$, so the transformed process $\widehat{X}$ has infinite lifetime and is conservative in the ergodic theory sense ("recurrent" in the sense of [8], page 48).*

PROOF. Our proof is a modification of the proof of Theorem 3.6 in [3]. For the reader's convenience, we spell out the details. Let $F_n^{(1)} = \{ x : \rho(x) \geq 1/n \}$. Then $\{ F_n^{(1)} : n \geq 1 \}$ is an $\mathcal{E}$-nest. By the probabilistic characterization of $\mathcal{E}$-nest, $\{ F_n^{(1)} \}_{n \geq 1}$ is an $\widehat{\mathcal{E}}$-nest (for $\widehat{X}$) as well. So there is an $\widehat{\mathcal{E}}$-nest $\{ F_n ; n \geq 1 \}$ of compact sets and a sequence of $g_n \in \widehat{\mathcal{F}}$ such that $F_n \subset F_n^{(1)}$ and $g_n = 1$ on $F_n$ for each $n \geq 1$. Again by the probabilistic characterization of $\mathcal{E}$-nest, $\{ F_n \}_{n \geq 1}$ is also an $\mathcal{E}$-nest (for $X$). If we let $\rho_{F_n}$ be the $\mathcal{E}_1$-orthogonal projection of $\rho$ onto $\mathcal{F}_{F_n} := \{ u \in \mathcal{F} : u = 0 \text{ q.e. on } F_n^c \}$, then $\rho_{F_n}$ converges to $\rho$ in $(\mathcal{F}, \mathcal{E}_1)$. Let $\rho_n = (0 \vee \rho_{F_n}) \wedge \rho$. Then $\rho_n$ converges to $\rho$ in $(\mathcal{F}, \mathcal{E}_1)$ as well (cf. Theorem 1.4.2(v) in [8]). Taking a subsequence if necessary, we may assume that $\rho_n$ converges to $\rho$, $\mathcal{E}$-q.e. on $E$. For $n \geq 1$, define $h_n = \rho_n / \rho$. Since $\rho \geq 1/n$ $\mathcal{E}$-q.e. on $F_n$ and $\rho_n = 0$ $\mathcal{E}$-q.e. on $F_n^c$, we have $h_n \in \mathcal{F}$, by the contraction property of $(\mathcal{E}, \mathcal{F})$. Note that $0 \leq h_n \leq 1$ and $h_n \to 1$ q.e. on $E$ as $n \to \infty$. By a calculation found in the proof of Lemma 6.3.3 of [8], it can be shown that

$$
(2.9) \quad \tfrac{1}{2} \int_E \rho(x)^2 \mu^c_{\langle h_n \rangle}(dx) + \int_{E \times E} (h_n(x) - h_n(y))^2 \rho(x) \rho(y) J(dx, dy) \to 0
$$

as $n \to \infty$.

Let $u$ be a bounded function in $\mathcal{F}$ with

$$\tfrac{1}{2}\int_E \rho(x)^2 \mu_{\langle u \rangle}^c(dx)$$
$$+ \int_{E \times E}(u(x) - u(y))^2 \rho(x)\rho(y)J(dx, dy) + \int u(x)^2 \rho(x)^2 m(dx) < \infty.$$

Fix $n \geq 1$ and define $f := uh_n$. Clearly $f$ is a bounded function in $\mathcal{F}$ satisfying the above inequality with $f$ in place of $u$. The process $f(X_t)$ admits the following Lyons–Zheng forward–backward martingale decomposition:

$$(2.10) \quad f(X_t) - f(X_0) = \tfrac{1}{2}(M_t^f - M_t^f \circ r_t), \qquad \mathbb{P}_m\text{-a.s. on } \{t < \zeta\},$$

where $M_t^f$ is the martingale part in Fukushima's decomposition of $f(X_t) - f(X_0)$. Recall that $d\widehat{\mathbb{P}}_x = Z_t^\rho d\mathbb{P}_x$ on $\mathcal{M}_t \cap \{t < \zeta_p\}$. Hence

$$K_t := M_t^f - \int_0^t \frac{1}{Z_{s-}^\rho} d\langle M^f, Z^\rho \rangle_s = M_t^f - \langle M^f, M \rangle_t, \qquad t < \zeta_p,$$

is a local martingale AF under $\widehat{\mathbb{P}} = (\widehat{\mathbb{P}}_x : x \in E)$ and

$$(2.11) \qquad [K]_t(\widehat{\mathbb{P}}) = [M^f]_t(\mathbb{P}), \qquad \widehat{\mathbb{P}}_m\text{-a.s. for } t < \zeta_p.$$

Here $[K](\widehat{\mathbb{P}})$ is the square bracket process for the martingale $K$ under the family $\widehat{\mathbb{P}}$, and $[M^f](\mathbb{P})$ is the square bracket for martingale $M^f$ under the family $\mathbb{P}$. We will use $\langle K \rangle(\widehat{\mathbb{P}})$ and $\langle M^f \rangle(\mathbb{P})$ to denote the dual predictable projections of $[K](\widehat{\mathbb{P}})$ and $[M^f](\mathbb{P})$ under the respective families $\widehat{\mathbb{P}}$ and $\mathbb{P}$. It follows from (2.11) that for $t < \zeta_p$,

$$\langle K \rangle_t(\widehat{\mathbb{P}}) = \langle M^f \rangle_t(\mathbb{P}) + \int_0^t \frac{1}{Z_{s-}^\rho} d\langle [M^f], Z^\rho \rangle_s(\mathbb{P}) = \langle M^f \rangle_t(\mathbb{P}) + \langle [M^f], M \rangle_t(\mathbb{P});$$

see, for example, Chapter 12 of [13]. By the quasi-left continuity of the processes $X$ and $\widehat{X}$, all of the "sharp bracket" processes involved are continuous, so we have, for $t \geq 0$,

$$\langle K \rangle_t(\widehat{\mathbb{P}}) = \langle M^f \rangle_t(\mathbb{P}) + \langle [M^f], M \rangle_t(\mathbb{P})$$

$$(2.12) \qquad = \langle M^f \rangle_t(\mathbb{P}) + \left( \sum_{0 < s \leq t}(f(X_s) - f(X_{s-}))^2 \left( \frac{\rho(X_s)}{\rho(X_{s-})} - 1 \right) \right)^p(\mathbb{P})$$

$$= \langle M^f \rangle_t(\mathbb{P}) + \int_0^t \int_{E_\Delta}(f(X_s) - f(y))^2 \left( \frac{\rho(y)}{\rho(X_s)} - 1 \right) N(X_s, dy)\, dH_s,$$

with the convention that $0/0 = 1$. Thus, by Theorem 2.5, the Revuz measure for the PCAF $\langle K \rangle_t(\widehat{\mathbb{P}})$ of $\widehat{X}$ is

$$\rho(x)^2 \mu_{\langle f \rangle}(dx)$$

$$(2.13) \qquad + 2\rho(x)^2 \int_{y \in E} (f(x) - f(y))^2 \Big(\frac{\rho(y)}{\rho(x)} - 1\Big) J(dx, dy) - f(x)^2 \rho(x)^2$$

$$= \rho(x)^2 \mu^c_{\langle f \rangle}(dx) + 2 \int_{y \in E} (f(x) - f(y))^2 \rho(x)\rho(y) J(dx, dy).$$

Now the CAF $\langle M^f, M \rangle$ is even, so by (2.10),

$$(2.14) \qquad f(X_t) - f(X_0) = \tfrac{1}{2}(K_t - K_t \circ r_t), \qquad \mathbb{P}_m\text{-a.s. on } \{t < \zeta\}.$$

Let $\nu = \rho^2 m$ and

$$\widehat{\mathbb{P}}_\nu(\cdot) = \int_E \widehat{\mathbb{P}}_x(\cdot)\nu(dx).$$

Applying Theorem 2.5 and noting that the time reversal operator $r_t$ also leaves the measure $\widehat{\mathbb{P}}_\nu$ invariant on $\mathcal{M}_t \cap \{t < \zeta\}$, we have by (2.13) and (2.14)

$$\lim_{t \to 0} \frac{1}{t} \widehat{\mathbb{E}}_\nu[(f(\widehat{X}_t) - f(\widehat{X}_0))^2; t < \zeta]$$

$$\leq \lim_{t \to 0} \Big( \frac{1}{2t} \widehat{\mathbb{E}}_\nu[(K_t)^2; t < \zeta] + \frac{1}{2t} \widehat{\mathbb{E}}_\nu[(K_t \circ r_t)^2; t < \zeta] \Big)$$

$$= \lim_{t \to 0} \frac{1}{t} \widehat{\mathbb{E}}_\nu[\langle K \rangle_t(\widehat{\mathbb{P}}); \ t < \zeta]$$

$$\leq \int_E \rho(x)^2 \mu^c_{\langle f \rangle}(dx) + 2 \int_{E \times E} (f(x) - f(y))^2 \rho(x)\rho(y) J(dx, dy)$$

$$< \infty.$$

Recall that $f = 0$ $m$-a.e. on $F_n^c$ and $g_n \in \widehat{\mathcal{F}}$ with $g_n = 1$ $m$-a.e. on $F_n$. Thus $f = fg_n$ and

$$\lim_{t \to 0} \frac{1}{t} \int_E (f(x) - \widehat{P}_t f(x)) f(x)\nu(dx)$$

$$= \lim_{t \to 0} \frac{1}{t} \Big( \frac{1}{2} \widehat{\mathbb{E}}_\nu[(f(\widehat{X}_t) - f(\widehat{X}_0))^2; \ t < \zeta] + \int_E f(x)^2 (1 - \widehat{P}_t 1)(x)\nu(dx) \Big)$$

$$\leq \limsup_{t \to 0} \frac{1}{2t} \widehat{\mathbb{E}}_\nu[(f(\widehat{X}_t) - f(\widehat{X}_0))^2; \ t < \zeta]$$

$$\quad + \limsup_{t \to 0} \frac{1}{t} \int_E (fg_n)(x)^2 (1 - \widehat{P}_t 1)(x)\nu(dx)$$

$$\leq \limsup_{t \to 0} \frac{1}{2t} \widehat{\mathbb{E}}_\nu[(f(\widehat{X}_t) - f(\widehat{X}_0))^2; \ t < \zeta]$$

$$\quad + \|f\|_\infty^2 \limsup_{t \to 0} \frac{1}{t} \int_E g_n(x)^2 (1 - \widehat{P}_t 1)(x)\nu(dx)$$

$$< \infty.$$

Therefore $f \in \widehat{\mathcal{F}}$, so $f$ admits a Fukushima decomposition:

$$f(\widehat{X}_t) - f(\widehat{X}_0) = \widehat{M}_t^f + \widehat{N}_t^f, \qquad \widehat{\mathbb{P}}_\nu\text{-a.s.,}$$

where $\widehat{M}_t^f$ is a $\widehat{\mathbb{P}}_x$-square-integrable martingale AF and $\widehat{N}_t^f$ is a CAF of zero energy; in particular, $\widehat{N}_t^f$ is a process of zero quadratic variation. On the other hand, $f(X_t)$ has a Fukushima decomposition under the family $\mathbb{P} = (\mathbb{P}_x : x \in E)$:

$$f(X_t) - f(X_0) = M_t^f + N_t^f.$$

By Girsanov's theorem the process $K_t = M_t^f - \langle M^f, M \rangle_t$ is a local martingale under $\widehat{\mathbb{P}}_x$ on $[\![0, \zeta_p[\![$, so by uniqueness we have

$$(2.15) \qquad \widehat{M}_t^f = K_t \qquad \text{for } t < \zeta_p.$$

To express $\widehat{\mathcal{E}}(f, f)$, we first calculate the killing measure $\widehat{\kappa}$ for the transformed process $\{\widehat{X}, \widehat{\mathbb{P}}_x, \ x \in E\}$. Now $\widehat{\kappa}$ is the Revuz measure of the PCAF $(\mathbb{1}_{\{t \geq \zeta_i\}})^p(\widehat{\mathbb{P}})$, the dual predictable projection of the increasing AF $t \to \mathbb{1}_{\{t \geq \zeta_i\}}$ under $\widehat{\mathbb{P}}$. By the same reasoning as for (2.12), for q.e. $x \in E$,

$$\begin{aligned}
(\mathbb{1}_{\{t \geq \zeta_i\}})^p(\widehat{\mathbb{P}}) &= (\mathbb{1}_{\{t \geq \zeta_i\}})^p(\mathbb{P}) + \langle \mathbb{1}_{\{t \geq \zeta_i\}}, M \rangle(\mathbb{P}) \\
&= (\mathbb{1}_{\{t \geq \zeta_i\}})^p(\mathbb{P}) - (\mathbb{1}_{\{t \geq \zeta_i\}})^p(\mathbb{P}) \\
&= 0.
\end{aligned}$$

Thus $\widehat{\kappa} = 0$. Now by (2.13) and (2.15),

$$\begin{aligned}
(2.16) \qquad \widehat{\mathcal{E}}(f, f) &= \lim_{t \to 0} \frac{1}{2t} \widehat{\mathbb{E}}_\nu[(f(\widehat{X}_t) - f(\widehat{X}_0))^2] \\
&= \lim_{t \to 0} \frac{1}{2t} \widehat{\mathbb{E}}_\nu[(\widehat{M}_t^f)^2] = \lim_{t \to 0} \frac{1}{2t} \widehat{\mathbb{E}}_\nu[\langle K \rangle_t] \\
&= \frac{1}{2} \int_E \rho(x)^2 \mu_{\langle f \rangle}^c(dx) + \int_{E \times E} (f(x) - f(y))^2 \rho(x)\rho(y) J(dx, dy).
\end{aligned}$$

Applying the above argument to $h_n$ (in place of $f$), we see that $h_n \in \widehat{\mathcal{F}}$ and, by (2.9), that $\widehat{\mathcal{E}}(h_n, h_n) \to 0$ as $n \to \infty$. Since $h_n \to 1$ q.e. on $E$, this implies that $1 \in \widehat{\mathcal{F}}_e \cap L^2(E, \rho^2 m) = \widehat{\mathcal{F}}$ and $\widehat{\mathcal{E}}(1, 1) = 0$ (see Theorem 1.5.2 of [8]). Consequently, $\widehat{X}$ is recurrent by Theorem 1.6.3 of [8]. This proves Theorem 2.6(b).

So far we have proved that $f = u h_n \in \widehat{\mathcal{F}}$ and that (2.16) holds for $f$. Note that

$$\begin{aligned}
\widehat{\mathcal{E}}(u h_n, u h_n) \leq{}& 2\|u\|_\infty^2 \widehat{\mathcal{E}}(h_n, h_n) + \int_E \rho(x)^2 \mu_{\langle u \rangle}^c(dx) \\
&+ 2 \int_{E \times E} (u(x) - u(y))^2 \rho(x)\rho(y) J(dx, dy),
\end{aligned}$$

which is uniformly bounded. As $|uh_n| \le |u|$, $uh_n \to u$, we see that $u$ can be approximated in $(\widehat{\mathcal{F}}, \widehat{\mathcal{E}}_1)$ by the Cesàro means of a subsequence of $\{uh_n\}_{n\ge 1}$. Hence $u$ is in $\widehat{\mathcal{F}}$. Repeating the computation for $f$ shows that (2.16) holds for $u$ as well. This proves Theorem 2.6(a). □

REMARK 2.7. Suppose that $\rho$ is in $\overset{\bullet}{\mathcal{F}}_{\mathrm{loc}}$ with $\rho > 0$ q.e. on $E$ and that $t \mapsto \sum_{0 < s \le t} |\rho(X_s) - \rho(X_{s-})|$ is locally $\mathbb{P}_x$-integrable for q.e. $x \in E$. This is the case if $\rho$ is bounded, for example. If we define

$$M_t^\rho = M_t^{\rho,c} + \left( \sum_{0 < s \le t} (\rho(X_s) - \rho(X_{s-})) \right)^{\sim},$$

where the superscript $^{\sim}$ indicates compensated sum, and use (2.1) and (2.2) to define $Z^\rho$, then Theorem 2.6(a) remains valid (with the same proof). We will not use this fact in the sequel.

We now identify the domain of the Dirichlet space for $\widehat{X}$.

THEOREM 2.8. *Under the condition of Theorem 2.6, the domain $\widehat{\mathcal{F}}$ of the Dirichlet form $(\widehat{\mathcal{E}}, \widehat{\mathcal{F}})$ for the Girsanov transformed process $\widehat{X}$ is the $\widehat{\mathcal{E}}_1$-completion of*

$$\left\{ f \in \mathcal{F} \colon \int \rho(x)^2 \mu_{\langle f \rangle}^c(dx) < \infty \right\} \cap \mathcal{L}^2(\rho \otimes \rho \cdot J) \cap L^2(\rho^2 m),$$

*where $\widehat{\mathcal{E}}_1 = \widehat{\mathcal{E}} + (\cdot, \cdot)_{L^2(E, \rho^2 m)}$.*

We first prepare a lemma.

LEMMA 2.9. *Define $\widehat{N}(x, dy) := \frac{\rho(y)}{\rho(x)} N(x, dy)$. Then $(N^{\widehat{X}}, H)$ is a Lévy system for $\widehat{X}$. Consequently, if $\widehat{J}$ denotes the jump measure of $\widehat{X}$, then $\widehat{J} = \rho \otimes \rho \cdot J$.*

PROOF. Let $K = (K_t)$ be a predictable process on $\Omega$ and let $f$ be a nonnegative Borel function on $E \times E$ vanishing on diagonal. By [24], (62.13) we find

$$\widehat{\mathbb{P}}_x \left[ \sum_{s \le t} K_s f(\widehat{X}_{s-}, \widehat{X}_s) \right] = \mathbb{P}_x \left[ \sum_{s \le t} K_s f(X_{s-}, X_s) Z_s^\rho; t < \zeta \right]$$

$$= \mathbb{P}_x \left[ \sum_{s \le t} K_s Z_{s-}^\rho f(X_{s-}, X_s) \frac{\rho(X_s)}{\rho(X_{s-})}; t < \zeta \right]$$

$$= \mathbb{P}_x \left[ \int_0^t K_s Z_s^\rho \rho(X_s)^{-1} N(\rho f)(X_s) \, dH_s; t < \zeta \right]$$

$$= \widehat{\mathbb{P}}_x \left[ \int_0^t K_s \widehat{N} f(\widehat{X}_s) \, dH_s \right].$$

The conclusion follows. □

PROOF OF THEOREM 2.8. Define $M^{\rho,o} = M^\rho + (\rho(X_{\zeta-}) \mathbb{1}_{\{\cdot \geq \zeta_i\}})^{\sim}$ and $M_t^o = \int_0^t [\rho(X_{s-})]^{-1} \, dM_s^{\rho,o}$. Clearly $M^{\rho,o}$ is a $\mathbb{P}_x$-square-integrable MAF of $X$ and $M^o$ is a locally $\mathbb{P}_x$-square-integrable MAF of $X$. Define

$$(2.17) \quad \widehat{M}_t := -M_t^o + \langle M^c \rangle_t(\mathbb{P}) + \sum_{0 < s \leq t} \frac{(\rho(\widehat{X}_s) - \rho(\widehat{X}_{s-}))^2}{\rho(\widehat{X}_s)\rho(\widehat{X}_{s-})} \mathbb{1}_{\{s < \zeta\}}.$$

Then $\widehat{M}$ is a local MAF of $\widehat{X}$. To see this, observe that, by Girsanov's theorem,

$$M_t^o - \langle M^o, M \rangle_t(\mathbb{P}) = M_t^o - \langle M^c \rangle_t(\mathbb{P}) - \left( \sum_{0 < s \leq \cdot} \left( \frac{\rho(X_s)}{\rho(X_{s-})} - 1 \right)^2 \mathbb{1}_{\{s < \zeta\}} \right)_t^p(\mathbb{P})$$

is a local martingale AF of $\widehat{X}$. An application of Lemma 2.9 [to compute the compensator of the sum on the right-hand side of (2.17)] now completes the proof of the claim.

Let $\widehat{Z}_t$ be the solution to $d\widehat{Z}_t = \widehat{Z}_{t-} \, d\widehat{M}_t$. Denote by $A^\kappa$ the PCAF of $X$ associated with the killing measure $\kappa$. Recall that $A^\kappa = (\mathbb{1}_{\{\cdot \geq \zeta_i\}})^p(\mathbb{P})$. Using the Doléans–Dade formula, one sees that

$$(2.18) \quad \frac{d\mathbb{P}_x}{d\widehat{\mathbb{P}}_x} \bigg|_{\mathcal{M}_t} = \frac{1}{Z_t^\rho} = \widehat{Z}_t \, e^{-A_t^\kappa}.$$

(cf. the last section in the proof of Theorem 3.7 in [3]).

Let $F_k = \{x \in E : \rho(x) \geq 1/k\}$, which is an $\mathcal{E}$-nest hence an $\widehat{\mathcal{E}}$-nest as well. Define $\rho_k = \rho - (\rho \wedge \frac{1}{k})$ and $h_k = \rho_k/\rho$. Clearly $0 \leq h_k \leq 1$, $\rho_k \in \mathcal{F}_{F_k}$, and $\rho_k \to \rho$ in $(\mathcal{F}, \mathcal{E}_1)$. Arguing as in the proof of Theorem 2.6, it can be shown that $h_k \in \mathcal{L}^2(\rho \otimes \rho \cdot J)$, and $h_k \in \widehat{\mathcal{F}}$ with $\widehat{\mathcal{E}}(h_k, h_k) \to 0$ as $k \to \infty$. Now for any $u \in b\widehat{\mathcal{F}}$, we claim $f := u h_k$ lies in $b\mathcal{F}_{F_k}$ and

$$\int \rho(x)^2 \mu_{\langle f \rangle}^c(dx) + \int (f(x) - f(y))^2 \rho(x)\rho(y) J(dx, dy) < \infty.$$

First note that

$$\int (f(x) - f(y))^2 J(dx, dy)$$

$$\leq \int_{F_k \times F_k} \bar{f}(x,y)^2 J(dx, dy) + \|u\|_\infty^2 \int_{(F_k \times F_k)^c} \bar{h}_k(x,y)^2 J(dx, dy)$$

$$\leq k^2 \int_{F_k \times F_k} \bar{f}(x,y)^2 \rho(x)\rho(y) J(dx,dy)$$

$$+ \|u\|_\infty^2 \int_{(F_k \times F_k)^c} \bar{h}_k(x,y)^2 J(dx,dy)$$

$$< \infty.$$

As $f \in b\widehat{\mathcal{F}}$, by Fukushima's decomposition,

$$f(\widehat{X}_t) - f(\widehat{X}_0) = \widehat{M}_t^f + \widehat{N}_t^f.$$

Define a family of measures $\mathbb{Q} = (\mathbb{Q}_x : x \in E)$ through $\frac{d\mathbb{Q}_x}{d\widehat{\mathbb{P}}_x}|_{\mathcal{M}_t} = \widehat{Z}_t$, where $\widehat{Z}$ is the Doléans–Dade exponential of $\widehat{M}$, and $\widehat{M}$ is as in (2.17). In this proof only we shall use $X^*$ to denote the coordinate process when referring to $\mathbb{Q}$. Recall that $\widehat{X}$ is a right Markov process with symmetry measure $\nu(dx) := \rho(x)^2 m(dx)$. It can be shown as in the proof of Lemma 2.4 that $(X^*, \mathbb{Q})$ is a right Markov process with symmetry measure $\rho(x)^{-2}\nu(dx) = m(dx)$. From Girsanov's theorem,

$$K := \widehat{M}^f - \langle \widehat{M}^f, \widehat{M} \rangle (\widehat{\mathbb{P}})$$

is a $\mathbb{Q}$-local martingale. Note that since $[K]_t(\mathbb{Q}) = [\widehat{M}^f]_t(\widehat{\mathbb{P}})$,

(2.19) $$\langle K \rangle_t(\mathbb{Q}) = \langle \widehat{M}^f \rangle_t(\widehat{\mathbb{P}}) + \langle [\widehat{M}^f], \widehat{M} \rangle_t(\widehat{\mathbb{P}}).$$

But

$$\langle [\widehat{M}^f], \widehat{M} \rangle_t(\widehat{\mathbb{P}}) = \left( \sum_{0 < s \leq t} (f(\widehat{X}_s) - f(\widehat{X}_{s-}))^2 \frac{\rho(\widehat{X}_{s-}) - \rho(\widehat{X}_s)}{\rho(\widehat{X}_s)} \right)^p (\widehat{\mathbb{P}}),$$

so its Revuz measure with respect to $(\widehat{X}, \widehat{\mathbb{P}})$ is, by Lemma 2.9,

$$\int_{y \in E} (f(x) - f(y))^2 \rho(x)(\rho(x) - \rho(y)) J(dx,dy).$$

Hence by applying Theorem 2.5 with $\widehat{X}$ and $X^*$ in the roles of $X$ and $\widehat{X}$, and with $\rho^{-1}$ in place of $\rho$, the Revuz measure of $\langle [\widehat{M}^f], \widehat{M} \rangle (\widehat{\mathbb{P}})$ viewed as a PCAF of $X^*$ is

$$\frac{1}{\rho(x)} \int_{y \in E} (f(x) - f(y))^2 (\rho(x) - \rho(y)) J(dx,dy).$$

Likewise, by (2.19), the Revuz measure $\mu_{\langle K \rangle}^*$ of $\langle K \rangle (\mathbb{Q})$ viewed as PCAF of $X^*$ is seen to be

(2.20) $$\mu_{\langle K \rangle}^*(dx) = \frac{1}{\rho(x)^2} \widehat{\mu}_{\langle f \rangle}^c(dx) + \int_{y \in E} (f(x) - f(y))^2 J(dx,dy).$$

As $\rho \geq 1/k$ on $F_k$ and $f = 0$ on $F_k^c$, hence $\int \rho(x)^{-2} \widehat{\mu}_{\langle f \rangle}^c(dx) < \infty$, and therefore $\mu_{\langle K \rangle}^*(E) < \infty$. But $f(X_t^*) - f(X_0^*) = (K_t - K_t \circ r_t)/2$, and we deduce, by reasoning similar to that used just below (2.14), that $f$ is in the Dirichlet space of $X^*$. In view of (2.18), we have $f = uh_k \in b\mathcal{F}_{F_k}$. As both $u$ and $h_k$ are in $b\widehat{\mathcal{F}}$, Lemma 2.9 yields

$$\int (f(x) - f(y))^2 \rho(x)\rho(y)J(dx,dy) < \infty.$$

By a calculation similar to that used in the proof of Lemma 2.9, the jump measure of $X^*$ is $J$. Likewise, by the argument appearing between (2.15) and (2.16), one sees that the killing measure of $X^*$ is the zero measure. If we use $\mu^*_{\langle f \rangle}$ and $\mu^{*,c}_{\langle f \rangle}$ to denote the energy measure of $f$ and its strong local part, in the context of $X^*$, then $\mu^*_{\langle f \rangle} = \mu^*_{\langle K \rangle}$. On the other hand,

$$\mu^*_{\langle f \rangle}(dx) = \widehat{\mu}^{*,c}_{\langle f \rangle}(dx) + \int_{y \in E} (f(x) - f(y))^2 J(dx,dy).$$

Hence from (2.20) we see that

$$\mu^{*,c}_{\langle f \rangle}(dx) = \frac{1}{\rho(x)^2}\widehat{\mu}^c_{\langle f \rangle}(dx).$$

As the Feynman–Kac transformation by the multiplicative functional $\exp(-A_t^\kappa)$ does not change the strongly local part of the energy measure, we have

$$\mu^c_{\langle f \rangle}(dx) = \frac{1}{\rho(x)^2}\widehat{\mu}^c_{\langle f \rangle}(dx),$$

and so

$$\int \rho(x)^2 \mu^c_{\langle f \rangle}(dx) + \int (f(x) - f(y))^2 \rho(x)\rho(y)J(dx,dy) < \infty.$$

Since $|uh_k| \leq |u|$, $uh_k \to u$ and

$$\widehat{\mathcal{E}}(uh_k, uh_k) \leq 2\|u\|_\infty^2 \widehat{\mathcal{E}}(h_k, h_k) + 2\widehat{\mathcal{E}}(u,u),$$

which is uniformly bounded, we see that $u$ can be approximated in $(\widehat{\mathcal{F}}, \widehat{\mathcal{E}}_1)$ by the Cesàro mean of a subsequence of $\{uh_k\}_{k \geq 1}$. Hence $u$ is in the $\widehat{\mathcal{E}}_1$-closure of

$$\left\{ f \in \mathcal{F} \colon \int \rho(x)^2 \mu^c_{\langle f \rangle}(dx) < \infty \right\} \cap \mathcal{L}^2(\rho \otimes \rho \cdot J) \cap L^2(\rho^2 m).$$

This proves the theorem.  $\square$

**3. Supermartingale multiplicative functional.** In this section we prove a representation theorem for a general class of supermartingale multiplicative functionals (MFs) of $X$. This result is a sharpening of results of Kunita ([16], Theorem 3.1) and Sharpe ([23], Theorem 7.1); see also [17], Section 6 and [15], Section 4. For a stopping time $T$, we will use $I(T)$ to denote the stochastic interval $[\![0, T[\![ \cup [\![T_i]\!]$, where $T_i$ is the totally inaccessible part of $T$. By a slight abuse of notation, we shall often write "$t \in I(T)$" to mean "$(t, \omega) \in I(T)$," where $\omega$ is the (suppressed, as usual) sample path.

THEOREM 3.1. *Let $Z$ be a supermartingale MF of $X$ such that $Z_0 \equiv 1$. Then there is a local martingale AF $M$, a PCAF $A$, and a Borel function $\varphi : E \times E_\Delta \to [-1, +\infty[$ such that*

$$(3.1) \quad \Delta M_t := M_t - M_{t-} = \varphi(X_{t-}, X_t) \qquad \forall t \in I(\zeta), \ \mathbb{P}_m\text{-}a.s.,$$

$$(3.2) \quad \begin{aligned} Z_t &= e^{M_t - (1/2)\langle M^c \rangle_t - A_t} \\ &\times \prod_{0 < s \le t} [1 + \varphi(X_{s-}, X_s)] e^{-\varphi(X_{s-}, X_s)} \qquad \forall t \in I(\zeta), \ \mathbb{P}_m\text{-}a.s. \end{aligned}$$

*The AF $M$ and the PCAF $A$ are determined by $Z$ up to $\mathbb{P}_m$-evanescence. In particular, $\varphi$ is uniquely determined by $Z$ modulo null sets of the measure $J^*(B) := J(B \cap (E \times E)) + \kappa(\pi_1(B \cap (E \times \{\Delta\})))$, where $\pi_1(x, y) := x$. Moreover,*

$$(3.3) \quad \int_0^t N(\mathbb{1}_{\{|\varphi| \le 1\}} \varphi^2 + \mathbb{1}_{\{\varphi > 1\}} \varphi)(X_s) \, dH_s + \int_0^t \varphi(X_s, \Delta) \, dA_s^\kappa < +\infty$$

*for all $t \in I(\zeta)$, $\mathbb{P}_m$-a.s. Finally,*

$$(3.4) \quad S := \inf\{t > 0 : Z_t = 0\} = \inf\{t > 0 : \varphi(X_{t-}, X_t) = -1\}$$

*$\mathbb{P}_m$-a.s. on $\{S < \zeta_p\}$.*

For a semimartingale $N$, let $\mathrm{Exp}(N)$ denote the unique solution $Y$ of

$$Y_t = 1 + \int_0^t Y_{s-} \, dN_s.$$

$\mathrm{Exp}(N)$ is called the stochastic exponential (in the sense of Doléans–Dade) of $N$. Formula (3.2) amounts to the statement that $Z = \mathrm{Exp}(M - A)$ at least on $I(\zeta)$. Before turning to the proof of Theorem 3.1 we prepare the way with a lemma.

LEMMA 3.2. *Let $B = (B_t)_{t \ge 0}$ be an AF of $X$. Then there is a Borel function $b : E \times E_\Delta \to \mathbb{R}$ with $b(x, x) = 0$ for all $x \in E$ such that*

$$(3.5) \quad \Delta B_t := B_t - B_{t-} = b(X_{t-}, X_t) \qquad \forall t \in ]0, \zeta_p[, \ \mathbb{P}_m\text{-}a.s.$$

*If $b'$ is another such function, then $J^*(b' \ne b) = 0$.*

PROOF. It follows from [11], (16.5) that there is a Borel function $b_0 : E \times E \to \mathbb{R}$ such that $\Delta B_t = b_0(X_{t-}, X_t)$ for all $t \in ]0, \zeta[$, $\mathbb{P}_m$-a.s. Fix $\epsilon > 0$ and define $T := \inf\{t > 0 : |\Delta B_t| \ge \epsilon\}$. Then $T$ is a *thin* terminal time; that is, $\mathbb{P}_x[T = 0] = 0$ for all $x \in E$. Consequently, $T_p$, the predictable part of $T$, is a thin predictable terminal time; by [11], (16.21), $\mathbb{P}_m(T_p < \zeta) = 0$. Owing to [24], Section 73, it follows that $\{t \in [0, \zeta[ : |\Delta B_t| > 0\} \subset \{t > 0 : X_{t-} \ne X_t\}$. Thus, modifying $b_0$ if necessary, we can arrange that $b_0(x, x) = 0$ for all

$x \in E$. Next, notice that $\Delta B_{\zeta_i}$ is measurable over the germ $\sigma$-field $\mathcal{F}^m_{[\zeta_i-, \zeta_i]}$. Because $X \equiv \Delta$ on $[\![\zeta_i, \infty[\![$, this germ $\sigma$-field is generated (modulo $\mathbb{P}_m$-null sets) by the random variables of the form $g(X)_{\zeta_i-}$ as $g$ varies over the bounded one-excessive functions of $X$; see [24], (24.32)(ii). But in the present context, natural AFs are continuous (by Corollary 3.17 in [6]), from which it follows that $g(X)_{\zeta_i-} = g(X_{\zeta_i-})$, $\mathbb{P}_m$-a.s. Therefore there is a Borel function $b_\Delta \colon E \to \mathbb{R}$ such that $\Delta B_{\zeta_i} = b_\Delta(X_{\zeta_i-})$, $\mathbb{P}_m$-a.s. on $\{\zeta_i < \infty\}$; see [11], (16.4). Defining $b(x, y) = b_0(x, y)\mathbb{1}_{E \times E}(x, y) + b_\Delta(x)\mathbb{1}_{\{\Delta\}}(y)$, we obtain the representation (3.5). The proof of the uniqueness assertion is left as an exercise to the reader. □

PROOF OF THEOREM 3.1. We begin with a discussion of the terminal time $S$ defined in (3.4); for related work see [10, 15, 26]. Clearly $S$ is a thin terminal time. Define the sequence $\{S^{(n)} : n \geq 1\}$ of iterates of $S$ by setting $S^{(1)} := S$ and $S^{(n+1)} := S^{(n)} + S \circ \theta_{S^{(n)}}$ for $n = 1, 2, \ldots$. [As a matter of convention, if $S^{(n)}(\omega) = +\infty$, then $S^{(k)}(\omega) = +\infty$ for all $k > n$.] Next define

$$C_t := \sum_n \mathbb{1}_{\{S^{(n)} \leq t\}}, \qquad t \geq 0$$

and

$$S^{(\infty)} := \uparrow \lim_n S^{(n)}.$$

Let $[\Delta]$ denote the sample path $\omega$ such that $X_t(\omega) = \Delta$ for all $t \geq 0$. With the convention $Z_t([\Delta]) \equiv 1$, we have $S([\Delta]) = +\infty$. Thus, if $S^{(n)}(\omega) = \zeta(\omega)$, then $S^{(k)}(\omega) = +\infty$ for all $k > n$. Consequently, $\{S^{(n)} : n \geq 1\}$ announces $S^{(\infty)}$ on $\{S^{(\infty)} < +\infty\}$. That is, $S^{(\infty)}$ is a thin predictable terminal time, hence $S^{(\infty)} \geq \zeta_p$, $\mathbb{P}_m$-a.s. by [11], (16.21). Since $\{C_t = \infty, t < \zeta\} = \{S^{(\infty)} \leq t < \zeta\}$, it follows that $C$ is finite on $I(\zeta)$, $\mathbb{P}_m$-a.s. Let us now apply Lemma 3.2 to $C$, taking into account the fact that $\Delta C$ takes values in $\{0, 1\}$ in the present situation. We find that there is a Borel set $\Lambda \subset E \times E_\Delta$, disjoint from the "diagonal" $\{(x, x) : x \in E\}$, such that

$$C_t = \sum_{0 < s \leq t} \mathbb{1}_\Lambda(X_{s-}, X_s)$$

provided $t \in I(\zeta)$, $\mathbb{P}_m$-a.s. In particular,

(3.6)     $$S = \inf\{t > 0 : (X_{t-}, X_t) \in \Lambda\}$$

$\mathbb{P}_m$-a.s. on $\{S < \zeta_p\}$. In fact, a little thought shows that the infimum in (3.6) is attained: $(X_{S-}, X_S) \in \Lambda$, $\mathbb{P}_m$-a.s. on $\{S < \zeta_p\}$.

From [14], Théorème 1 and the quasi-left-continuity of the filtration $(\mathcal{M}_t)$, we know that $Z^{-1}$ is locally bounded on $I(S)$ in the sense that there is an increasing sequence $(S_n)$ of stopping times with $I(S) = \bigcup_{n=1}^{\infty} [\![0, S_n]\!]$ such that $Z_-^{-1} \leq n$ on $[\![0, S_n]\!]$ for all $n$, $\mathbb{P}_m$-a.s. Let $Z = N - V$ be the Doob–Meyer decomposition of $Z$ into local martingale and predictable increasing components $N$ and $V$, respectively. Define $M_t^\circ := \int_0^t Z_{s-}^{-1} \, dN_s$ and $A_t^\circ := \int_0^t Z_{s-}^{-1} \, dV_s$. Both of these integrals are well defined on $I(S)$ because of the local boundedness of $Z_-^{-1}$ just noted. It is not hard to check that $M^\circ$ is an AF of $(X, S)$ (the process $X$ killed at time $S$, defined to equal $X$ before $S$ and $\Delta$ at and after $S$) and a local martingale on $I(S)$, and that $A^\circ$ is a predictable increasing AF of $(X, S)$. Of course, $A^\circ$ is continuous except perhaps for a jump at $S_p$ on $\{S_p = \zeta_p < \infty\}$. (By the discussion in the preceding section, $S_{\{S < \zeta_p\}}$ is totally inaccessible.)

We now extend $M^\circ$ and $A^\circ$ to AFs of $X$. Let us begin with $A^\circ$. By an extension procedure detailed in Theorem (4.8) of [11], there is a diffuse homogeneous random measure $\alpha$ of $X$ such that $A_t^\circ = \alpha([0, t])$ for all $t < S$, $\mathbb{P}_m$-a.s. We will show that $\alpha([0, t]) < +\infty$ for all $t < \zeta$, $\mathbb{P}_m$-a.s., and then use the recipe $A_t := \alpha([0, t \wedge \zeta])$, $t \geq 0$, as the desired extension of $A^\circ$. Define $R := \inf\{t : \alpha([0, t]) = +\infty\}$. Clearly $R$ is a terminal time, and $R \geq S$, so that $R$ is thin and hence *exact*, in the sense that $t + R \circ \theta_t$ decreases to $R$ as $t$ decreases to 0. Now define $R_n := \inf\{t : \alpha([0, t]) \geq n\}$, and notice that $R_n \leq R$. The strong Markov property (applied at the stopping time $R$) shows that the event $\{R_n = R < \zeta\}$ differs from the event $\{R_n = R < \zeta, S \circ \theta_R > 0\}$ by a $\mathbb{P}_m$-null set. Suppose that $\omega$ is a point of $\{R_n = R < \zeta, S \circ \theta_R > 0\}$. Then $\alpha(\omega, [0, t]) < n$ for $0 \leq t < R(\omega)$ but $\alpha(\omega, [0, t]) = +\infty$ for $t > R(\omega)$. In particular, if $t > R(\omega)$, then

$$+\infty = \alpha(\omega, [0, t]) = \alpha(\omega, [0, R(\omega)]) + \alpha(\theta_R \omega, (0, t - R(\omega)]).$$

The furthest right term above is finite when $t$ is sufficiently close to $R(\omega)$ since $S(\theta_R \omega) > 0$. Therefore $\alpha(\omega, [0, R(\omega)]) = +\infty$. On the other hand, $\alpha(\omega, [0, R(\omega))) = \alpha(\omega, [0, R_n(\omega))) \leq n$. But $\alpha(\omega, \{R(\omega)\}) = 0$ because $\alpha$ is diffuse. It follows that $\mathbb{P}_m(R_n = R < \zeta) = 0$, so that $R' := R_{\{R < \zeta\}}$ is a thin *predictable* terminal time. Thus $R' \geq \zeta$, which forces $R \geq \zeta$ as well. This shows that $t \mapsto \alpha([0, t])$ defines a PCAF of $X$.

Turning to $M^\circ$, notice that $\Delta M_t^\circ = \Delta Z_t / Z_{t-} = Z_t / Z_{t-} - 1 \geq -1$. Define $B_t := \sum_{0 < s \leq t} \Delta M_s^\circ \mathbb{1}_{\{\Delta M_s^\circ > 1\}}$, let $B_t^p$ denote the dual predictable projection of $B$, and now define $M_t^{\circ,2} := B_t - B_t^p$ and $M_t^{\circ,1} := M_t^\circ - M_t^{\circ,2}$. Both $M^{\circ,1}$ and $M^{\circ,2}$ are local martingales on $I(S)$ and AFs of $(X, S)$. By the argument of the last section, the CAF [of $(X, S)$] $B^p$ extends to a CAF $\bar{B}$ of $X$. Moreover, by (3.5) there is a Borel function $\varphi : E \times E_\Delta \to [-1, +\infty]$ such that $\Delta M^\circ = \varphi(X_-, X)$ on $I(S)$. Then

$$M_t^{*,2} := \sum_{0 < s \leq t} \varphi(X_{s-}, X_s) \mathbb{1}_{\{\varphi > 1\}}(X_{s-}, X_s) - \bar{B}_t$$

defines an extension of $M^{\circ,2}$ to a local martingale AF of $X$. Next, given a locally square-integrable martingale AF $N$, consider the covariation process $\Phi^\circ(N) := \langle M^{\circ,1}, N \rangle$, viewed as a CAF of $(X, S)$. As before, this CAF admits a unique extension $\Phi(N)$ to a CAF of $X$. In addition, we have the Kunita–Watanabe estimate

$$\langle M^{\circ,1}, N \rangle_t^2 \leq \langle M^{\circ,1} \rangle_t \cdot \langle N \rangle_t \qquad \forall t \in [0, S[,$$

and it is easy to check that $[\Phi(N)_t]^2 \leq D_t \cdot \langle N \rangle_t$, where $D$ is the extension of the CAF $\langle M^{\circ,1} \rangle$. A result of Kunita (Proposition 2.4 in [15]) now tells us that there is a local martingale AF $M^{*,1}$ such that $\Phi(N) \equiv \langle M^{*,1}, N \rangle$ for all $N$. Of course, $M^{*,1} \equiv M^{\circ,1}$ on $I(S)$. The local martingale AF $M := M^{*,1} + M^{*,2}$ is the desired extension of $M^\circ$. Notice that $\Delta M \equiv \varphi(X_-, X)$ on $I(\zeta)$.  $\square$

## 4. Absolute continuity and Dirichlet forms.
Let $Y = (\Omega, \mathcal{M}, \mathcal{M}_t, Y_t, \mathbb{Q}_x, x \in E)$ be another symmetric Markov process with symmetry measure $\nu$, which is realized on the same (canonical) path space $\Omega$ as $(X, \mathbb{P}_x, x \in E)$. Here $Y_t(\omega) = X_t(\omega)$ but we use $Y$ for emphasis when referring to $Y$. As with the process $X$, we assume that $Y$ is a Borel right process.

We note that Lemmas 2.5, 3.4 and 3.9 in [5] are valid in the setting of symmetric Borel right processes. The first result of this section is the analogue of [5], Theorem 3.2.

THEOREM 4.1. *Assume* $\mathbb{Q}_\nu \ll_{\mathrm{loc}} \mathbb{P}_m$, *in that* $\mathbb{Q}_\nu|_{\mathcal{M}_t \cap \{t < \zeta\}}$ *is absolutely continuous with respect to* $\mathbb{P}_m|_{\mathcal{M}_t \cap \{t < \zeta\}}$ *for each* $t > 0$. *Then* $\nu \ll m$ *and there is a* $(Y, \nu)$-*inessential Borel set* $N \subset E$ *which is* $X$-*finely closed, and a version* $\rho^2$ *of the Radon–Nikodym derivative* $d\nu/dm$ *such that* $0 < \rho(x) < \infty$ *for all* $x \in E \setminus N$ *and:*

(a) $t \to \rho(X_t)$ *is right-continuous on* $[0, T_N[$ *with left limits on* $]0, T_N \wedge \zeta[$, $P_x$-*a.s. for all* $x \in E \setminus N$; *in particular,* $\rho|_{E \setminus N}$ *is* $(X, T_N)$-*finely continuous;*

(b) $\log \rho \in \overset{\bullet}{\mathcal{F}}_{\mathrm{loc}}(X, T_N)$ *and there exists a local martingale AF* $M_t$ *satisfying*

$$\log \rho(X_t) - \log \rho(X_0) = (M_t - M_t \circ r_t)/2, \qquad \mathbb{P}_m\text{-a.s. on } \{t < T_N \wedge \zeta\}.$$
(4.1)

PROOF. The existence of $N$ follows from the proofs of Lemmas 3.4 and 3.9 in [5]. By Kunita ([16], Theorem 5.1), there is a supermartingale multiplicative functional $Z_t$ of $X$ satisfying

$$(4.2) \qquad \frac{d\mathbb{Q}_\nu}{d\mathbb{P}_m}\bigg|_{\mathcal{M}_t \cap \{t < \zeta\}} = Z_t.$$

We set $Z_t = 0$ for $t \geq \zeta$. Since $X$ is symmetric under the measure $\mathbb{P}_m$, we have

$$\rho^2(X_0) Z_t = \rho^2(X_t) Z_t \circ r_t, \qquad \mathbb{P}_m\text{-a.s. on } \{t < \zeta\}.$$

This implies that

$$(4.3) \qquad 2\log\rho(X_t) - 2\log\rho(X_0)$$
$$= \log Z_t - \log Z_t \circ r_t, \qquad \mathbb{P}_m\text{-a.s. on } \{t < T_N \wedge \zeta\}.$$

But, by Theorem 3.1,

$$\log Z_t - \log Z_0 = M_t - \tfrac{1}{2}\langle M^c\rangle_t - A_t + H_t, \qquad 0 \le t < \zeta,$$

where

$$H_t := \sum_{0 < s \le t}\left(\log(1 + \varphi(X_{s-}, X_s)) - \varphi(X_{s-}, X_s)\right),$$

which is absolutely convergent in view of (3.3). Clearly $H$ is quasi-left-continuous since $X$ is so; thus its dual predictable projection $H^p$ is a CAF. Define

$$M_t^* := \tfrac{1}{2}(M_t + H_t - H_t^p), \qquad A_t^* := \tfrac{1}{2}(H_t^p - \tfrac{1}{2}\langle M^c\rangle_t - A_t).$$

Then $M_t^*$ is a local martingale AF, and $A_t^*$ is a CAF of finite variation with

$$\log Z_t - \log Z_0 = 2M_t^* + 2A_t^*.$$

Note that $Z_0 = \rho^2(X_0)$ and $A_t^*$ is even. So on $\{t < T_N \wedge \zeta\}$,

$$\log Z_t - \log Z_t \circ r_t + 2\log\rho(X_t) - 2\log\rho(X_0) = 2M_t^* - 2M_t^* \circ r_t.$$

Hence by (4.3)

$$\log\rho(X_t) - \log\rho(X_0) = (M_t^* - M_t^* \circ r_t)/2.$$

Theorem 4.2 now implies that $\log\rho \in \overset{\bullet}{\mathcal{F}}_{\mathrm{loc}}(X, T_N)$. Therefore $\rho$ has a quasi-continuous version with respect to subprocess $(X, T_N)$. $\square$

**THEOREM 4.2.** *Let $f$ be a quasi-continuous Borel function on $E$, and suppose that there is a local martingale AF $M$ such that*

$$f(X_t) - f(X_0) = (M_t - M_t \circ r_t)/2, \qquad \mathbb{P}_m\text{-a.s. for each } t \in [0, \zeta[.$$

*Then $f \in \overset{\bullet}{\mathcal{F}}_{\mathrm{loc}}$.*

PROOF. We first assume that the jumps of $M$ are bounded, and then we show how to reduce to this special case.

If there is a constant $C$ such that $|\Delta M_t(\omega)| \le C$, then it is easy to check that $[M]$ is locally integrable, so $\langle M\rangle$ exists and is a PCAF of $X$. Things being so, the argument in [5], Lemma 3.15, can be used to reach the desired conclusion.

In general, define $T := \inf\{t > 0 : |\Delta M_t| > 1\}$. Then $T$ is a thin terminal time, and the subprocess $(X, T)$ is $m$-symmetric with state space $E$. Theorem 4.1 in [27] provides a precise description of the Dirichlet form of $(X, T)$, telling us, in particular, that the Dirichlet space $\mathcal{F}(X, T)$ of $(X, T)$ is a subspace of $\mathcal{F}$. Evidently, $X$ and $(X, T)$ have the same fine topologies (modulo $X$-exceptional sets). Clearly, if $N$ is $X$-exceptional, then it is $(X, T)$-exceptional. Conversely, if $N$ is $(X, T)$-exceptional, then it is $X$-exceptional. To see this let $\{T_n\}$ denote the sequence of iterates of $T$. By assumption, $X$ does not encounter $N$ during any of the open intervals $]T_n, T_{n+1}[$. Thus $X$ visits $N$ at most countably often, $\mathbb{P}_m$-a.s. That is, $N$ is $m$-semipolar, hence exceptional since $X$ is symmetric. Using quasi-left-continuity, one now checks that any increasing sequence $\{G_n\}$ of finely open sets is an $X$-nest if and only if it is an $(X, T)$-nest.

Now for $f$ under the assumptions of the theorem, by modifying $M$ at time $T$ we can produce a local martingale AF of $(X, T)$, call it $M^*$, with jumps bounded by 1, such that

$$f(X_t) - f(X_0) = (M_t^* - M_t^* \circ r_t)/2, \qquad \mathbb{P}_m\text{-a.s. for each } t \in [0, T[.$$

By the first section, $f \in \overset{\bullet}{\mathcal{F}}(X, T)_{\mathrm{loc}}$. The preceding section tells us that $\overset{\bullet}{\mathcal{F}}(X, T)_{\mathrm{loc}} \subset \overset{\bullet}{\mathcal{F}}_{\mathrm{loc}}$ since $\mathcal{F}(X, T) \subset \mathcal{F}$. This completes the proof. $\quad\square$

For $t > 0$, we say that two sample paths $\omega$ and $\omega'$ are *pre-$t$-equivalent* provided $\omega(s) = \omega'(s)$ for all $s \in [0, t[$. Observe that if $A = (A_t)$ is a finite CAF of $X$ and if $\omega$ and $\omega'$ are pre-$t$-equivalent, then

$$A_s(\omega) = A_s(\omega') \qquad \text{for all } 0 \le s \le t.$$

It is easy to check that $r_t \theta_s \omega$ is pre-$t$-equivalent to $r_{t+s}\omega$ and that $\theta_t r_{t+s}\omega$ is pre-$s$-equivalent to $r_s \omega$. This will be used repeatedly in the proof of next theorem. Define

$$\widehat{A}_t = A_t \circ r_t \qquad \text{on } \{t < \zeta\}.$$

Following [26], we have the following result.

THEOREM 4.3.  $\widehat{A} = (\widehat{A}_t : 0 \le t < \zeta)$ *is a CAF of* $X$.

PROOF.  First we need to show that $\widehat{A}_t$ is an AF. On $\{t + s < \zeta\}$,

$$\widehat{A}_{t+s} = A_{t+s} \circ r_{t+s} = (A_t + A_s \circ \theta_t) \circ r_{t+s} = A_t \circ r_{t+s} + A_s \circ \theta_t \circ r_{t+s}$$
$$= A_t \circ r_t \circ \theta_s + A_s \circ r_s = \widehat{A}_s + \widehat{A}_t \circ \theta_s.$$

Note that on $\{t < \zeta\}$, for $0 < u < t$,

$$\widehat{A}_t - \widehat{A}_{t-u} = \widehat{A}_u \circ \theta_{t-u} = A_u \circ r_u \circ \theta_{t-u} = A_u \circ r_t.$$

Hence

$$\lim_{u \downarrow 0} (\widehat{A}_{t-u} - \widehat{A}_t) = -\lim_{u \downarrow 0} A_u \circ r_t = 0.$$

This shows that $\widehat{A}$ is left-continuous.

Let us now prove the right-continuity. Note that on $\{t + u < \zeta\}$,

$$\widehat{A}_{t+u} - \widehat{A}_t = \widehat{A}_u \circ \theta_t = A_u \circ r_u \circ \theta_t = A_u \circ r_{t+u},$$

so it suffices to show that $\lim_{u \downarrow 0} A_u \circ r_{t+u} = 0$. For any $s > u > 0$, since $\theta_{s-u} \circ r_{t+s} \omega$ is pre-$(t+u)$-equivalent to $r_{t+u} \omega$, we have

$$(A_s - A_{s-u}) \circ r_{t+s} = A_u \circ \theta_{s-u} \circ r_{t+s} = A_u \circ r_{t+u}.$$

Thus

$$\lim_{u \downarrow 0} A_u \circ r_{t+u} = \lim_{u \downarrow 0} (A_s - A_{s-u}) \circ r_{t+s} = (A_s - A_{s-}) \circ r_{t+s} = 0.$$

This proves the theorem. $\square$

For simplicity, from now on we will assume that $\mathbb{Q}_x \ll_{\mathrm{loc}} \mathbb{P}_x$ for all $x$ and that $Z > 0$ on $[\![0, \zeta[\![$. It is easy to reduce to this case by killing $X$ at a terminal time and removing a $Y$-exceptional set $N$ from $E$.

Defining $\ell := \log \rho$, we have $\ell \in \overset{\bullet}{\mathcal{F}}_{\mathrm{loc}}$ by Theorem 4.1. Recall that the density process $Z$ in (4.2) is a nonnegative supermartingale MF of $X$, which is strictly positive on $[\![0, \zeta[\![$. Hence by Theorem 3.1, $Z = \mathrm{Exp}(M - A)$, where $M$ is a local martingale MF and $A$ is a PCAF of $X$. Let $M^c$ and $M^d$ be the continuous and purely discontinuous components of $M$, and let $\varphi$ be the Borel function: $E \times E_\Delta \to [-1, +\infty[$ with $\varphi(x, x) = 0$ for all $x \in E$ such that

$$\Delta M_t^d = \Delta M_t = \varphi(X_{t-}, X_t), \qquad \mathbb{P}_m\text{-a.s.}$$

In particular, $M^d$ is the compensated local martingale corresponding to $\sum_{0 < s \leq .} \varphi(X_{s-}, X_s)$. We now deduce from the identity (4.1) that

$$(4.4) \quad \ell(X_t) - \ell(X_0) = \frac{1}{2}(M_t^c - M_t^c \circ r_t) + \frac{1}{2} \sum_{0 < s \leq t} \log\left(\frac{1 + \varphi(X_{s-}, X_s)}{1 + \varphi(X_s, X_{s-})}\right)$$

$\mathbb{P}_m$-a.s. on $\{t < \zeta\}$ for every $t > 0$. The infinite series in (4.4) is to be understood in the following sense:

$$
\begin{aligned}
(4.5) \quad & \sum_{0 < s \leq t} \log\left(\frac{1 + \varphi(X_{s-}, X_s)}{1 + \varphi(X_s, X_{s-})}\right) \\
&= \sum_{0 < s \leq t} \left(\log(1 + \varphi(X_{s-}, X_s)) - \varphi(X_{s-}, X_s)\right) \\
&\quad - \sum_{0 < s \leq t} \left(\log(1 + \varphi(X_s, X_{s-})) - \varphi(X_s, X_{s-})\right),
\end{aligned}
$$

both sums on the right being absolutely convergent, $\mathbb{P}_m$-a.s. on $\{t < \zeta\}$. By Theorem 4.3, (4.4) in fact holds for all $t \in [0, \zeta[$, $\mathbb{P}_m$-a.s.

As $\ell \in \overset{\bullet}{\mathcal{F}}_{\mathrm{loc}}$ by Theorem 4.1, it follows from Lemma 2.2 that

$$
\begin{aligned}
(4.6) \quad \ell(X_t) - \ell(X_0) &= \tfrac{1}{2}(M_t^{\ell,c} - M_t^{\ell,c} \circ r_t) \\
&\quad + \lim_{\varepsilon \downarrow 0} \sum_{0 < s \le t} (\ell(X_s) - \ell(X_{s-})) \mathbb{1}_{\{|\ell(X_s) - \ell(X_{s-})| > \varepsilon\}}
\end{aligned}
$$

for all $t \in [0, \zeta[$, $\mathbb{P}_m$-a.s. Identities (4.4)–(4.6) yield $\mathbb{P}_m$-a.s.

$$
M_t^c - M_t^c \circ r_t = M_t^{\ell,c} - M_t^{\ell,c} \circ r_t \qquad \text{for all } t \in [0, \zeta[.
$$

Using the fact that an even martingale CAF must vanish [5], (3.25), we deduce from the above that

$$
M^c = M^{\ell,c} = \rho^{-1}(X_-) \bullet M^{\rho,c},
$$

where $\rho^{-1}(X_-) \bullet M^{\rho,c}$ is an Itô integral, and the last equality follows from [8], Theorem 5.6.2. It follows from Lemma 3.2 and (4.5) that

$$
\ell(x) - \ell(y) = \frac{1}{2} \log\left( \frac{1 + \varphi(x,y)}{1 + \varphi(y,x)} \right) \qquad J^*\text{-a.e. on } E \times E,
$$

and so

$$
\frac{\rho(y)^2}{\rho(x)^2} = \frac{1 + \varphi(x,y)}{1 + \varphi(y,x)}.
$$

Thus, the function $\gamma$ defined by

$$
(4.7) \qquad \gamma(x,y) := \rho(x)^2 [1 + \varphi(x,y)]
$$

is symmetric, $J^*$-a.e.

Recall that by Theorem 3.1 the density $Z_t$ in (4.2) can be written as

$$
Z = \mathrm{Exp}(M - A) = \mathrm{Exp}(M)\,\mathrm{Exp}(-A),
$$

where $M$ is a local martingale MF and $A$ is a PCAF of $X$. As $\Delta Z_t = Z_{t-} \cdot \Delta M_t$, we have $Z_t = Z_{t-}(1 + \Delta M_t) = Z_{t-}(1 + \varphi(X_{t-}, X_t))$. By an argument used in the proof of Lemma 2.9, one sees that $(N^Y, H)$ is a Lévy system of $Y$, where $N^Y(x, dy) := (1 + \varphi(x,y))N(x, dy)$. Hence the jump measure of $Y$ is

$$
(4.8) \qquad J^Y(dx, dy) = \rho(x)^2(1 + \varphi(x,y))J(dx, dy) = \gamma(x,y)J(dx, dy)
$$

(cf. [5], Lemma 4.4).

To find the killing measure $\kappa^Y$ of $Y$, note that $\mathbb{Q} = (\mathbb{Q}_x : x \in E)$ and $\mathbb{P} = (\mathbb{P}_x : x \in E)$ are related by first making a Girsanov transformation using $\mathrm{Exp}(M)$ and then a Feynman–Kac transformation using $\mathrm{Exp}(-A)$. Hence $\kappa^Y$ is the Revuz measure for the PCAF

$$(\mathbb{1}_{\{\cdot \geq \zeta_i\}})^p(\mathbb{Q}) = (\mathbb{1}_{\{\cdot \geq \zeta_i\}})^p(\mathbb{P}) + \langle \mathbb{1}_{\{\cdot \geq \zeta_i\}}, M \rangle(\mathbb{P}) + A$$
$$= (\mathbb{1}_{\{\cdot \geq \zeta_i\}})^p(\mathbb{P}) + (\Delta M_{\zeta_i} \mathbb{1}_{\{\cdot \geq \zeta_i\}})^p(\mathbb{P}) + A$$
$$= (\mathbb{1}_{\{\cdot \geq \zeta_i\}})^p(\mathbb{P}) + (\varphi(X_{\zeta_i-}, \Delta) \mathbb{1}_{\{\cdot \geq \zeta_i\}})^p(\mathbb{P}) + A.$$

So

$$(4.9) \qquad \kappa^Y(dx) = \rho(x)^2(1 + \varphi(x, \Delta))\kappa(dx) + \rho(x)^2 \mu_A(dx).$$

The above discussion proves half of the theorem to follow. Let $(\mathcal{E}^Y, \mathcal{F}^Y)$ be the Dirichlet form for process $(Y, \mathbb{Q})$. From its probabilistic characterization, $\{F_k\}_{k\geq 1}$ is an $\mathcal{E}^Y$-nest if and only if it is an $\mathcal{E}$-nest. So we can choose $\mathcal{E}^Y$-nest $\{F_k\}_{k\geq 1}$ of compact sets (see, e.g., VI.3 of [19]) such that $1/k \leq \rho \leq k$ q.e. on $F_k$.

THEOREM 4.4. *Let $(\mathcal{E}^Y, \mathcal{F}^Y)$ and $\{F_k\}_{k\geq 1}$ be as above. Then*

$$(4.10) \qquad \bigcup_{k=1}^{\infty} \mathcal{F}_{F_k} \cap \mathcal{L}^2(\gamma \cdot J) \cap L^2(\kappa^Y) \subset \mathcal{F}^Y$$

*and for $f$ in the set on the left-hand side of* (4.10),

$$(4.11) \qquad \begin{aligned} \mathcal{E}^Y(f, f) = {} & \tfrac{1}{2} \int \rho(x)^2 \mu_{\langle f \rangle}^c(dx) \\ & + \int (f(x) - f(y))^2 J^Y(dx, dy) + \int f(x)^2 \kappa^Y(dx), \end{aligned}$$

*with $J^Y$ and $\kappa^Y$ given by* (4.8) *and* (4.9).

PROOF. To calculate the continuous part of the energy measure, we proceed as in Section 2, using the method of forward–backward martingale decompositions together with martingale theory. Define $\widetilde{\mathbb{Q}} = (\widetilde{\mathbb{Q}}_x : x \in E)$ by $d\widetilde{\mathbb{Q}}_x/d\mathbb{P}_x|_{\mathcal{M}_t \cap \{t < \zeta_p\}} = \mathrm{Exp}(M_t)$. Then $(X, \widetilde{\mathbb{Q}})$ is $\nu$-symmetric. For $f \in b\mathcal{F}_{F_k} \cap L^2(\kappa^Y) \cap \mathcal{L}^2(\gamma \cdot J)$, by the Lyons–Zheng forward–backward martingale decomposition,

$$f(X_t) - f(X_0) = \tfrac{1}{2}(M_t^f - M_t^f \circ r_t), \qquad \mathbb{P}_m\text{-a.s. on } \{t < \zeta\},$$

where $M_t^f$ is the martingale part in Fukushima's decomposition of $f(X_t) - f(X_0)$. Hence

$$K_t := M_t^f - \langle M^f, M \rangle_t, \qquad t < \zeta_p,$$

is a local martingale AF under $\widetilde{\mathbb{Q}}$ and

$$[K]_t(\widetilde{\mathbb{Q}}) = [M^f]_t(\mathbb{P}), \qquad \widehat{\mathbb{P}}_x\text{-a.s. for } t < \zeta_p.$$

Hence for $t < \zeta_p$,

$$\langle K \rangle_t(\widetilde{\mathbb{Q}}) = \langle M^f \rangle_t(\mathbb{P}) + \langle M^f, M \rangle_t(\mathbb{P})$$

$$= \langle M^f \rangle_t(\mathbb{P}) + \int_0^t \int_{E_\Delta} (f(X_s) - f(y))^2 \varphi(X_s, y) N(X_s, dy)\, dH_s.$$

So the Revuz measure of $\langle K \rangle$ with respect to $(X, \widetilde{\mathbb{Q}})$ is

$$\rho(x)^2 \mu_{\langle f \rangle}(dx) + 2\rho(x)^2 \int_{y \in E} (f(x) - f(y))^2 \varphi(x, y) J(dx, dy)$$

$$+ \rho(x)^2 \varphi(x, \Delta) f(x)^2 \kappa(dx)$$

$$= \rho(x)^2 \mu_{\langle f \rangle}^c(dx) + 2 \int_{y \in E} (f(x) - f(y))^2 J^Y(dx, dy)$$

$$+ \rho(x)^2 (1 + \varphi(x, \Delta)) f(x)^2 \kappa(dx).$$

Let $(\widehat{\mathcal{E}}, \widehat{\mathcal{F}})$ [resp. $(\mathcal{E}^Y, \mathcal{F}^Y)$] be the Dirichlet form for the process $(X, \widetilde{\mathbb{Q}})$ [resp. $(Y, \mathbb{Q})$]. As $(Y, \mathbb{Q})$ is obtained from $(X, \widetilde{\mathbb{Q}})$ through Feynman–Kac transform by $\mathrm{Exp}(-A)$, $\mathcal{F}^Y = \widehat{\mathcal{F}} \cap L^2(\phi^2 \mu_A)$ and $\mathcal{E}^Y(f, f) = \widehat{\mathcal{E}}(f, f) + \int f(x)^2 \phi(x)^2 \times \mu_A(dx)$. Because $f = 0$ q.e. on $F_k^c$ for some $k \geq 1$, an argument similar to that used in the proof of Theorem 2.6 [between (2.14) and (2.15)] shows that $f \in \widehat{\mathcal{F}}$. Applying Feynman–Kac, one has $f \in \mathcal{F}^Y$ with

$$\mathcal{E}^Y(f, f) = \tfrac{1}{2} \int \rho(x)^2 \mu_{\langle f \rangle}^c(dx) + \int (f(x) - f(y))^2 J^Y(dx, dy)$$

$$+ \int f(x)^2 \kappa^Y(dx),$$

where $J^Y$ and $\kappa^Y$ are given by (4.8) and (4.9).   □

In the remainder of this section, we will focus on the special case in which the supermartingale $Z_t$ in (4.2) is of pure jump type. That is, in the expression (3.2) for $Z$, we assume that $M^c = 0$, $A = 0$, $\varphi$ is symmetric on $E \times E$ with $\varphi > -1$ on $E \times E$ and $\varphi(x, \Delta) = 0$, so $Z$ is strictly positive on $[\![0, \zeta_p[\![$ and $t \mapsto Z_t$ is continuous at $\zeta_i$. As a consequence of Theorem 3.1, the integrability condition (3.3) holds.

COROLLARY 4.5.  *In the setting of Theorem 4.4, suppose there are real constants $c_1$ and $c_2$ such that $-1 < c_1 \leq \varphi(x, y) \leq c_2$ for all $x, y \in E$. Then $\mathcal{F}^Y = \mathcal{F}$ and (4.11) holds for all $f \in \mathcal{F}$.*

PROOF.  By Theorem 4.4, $\mathcal{F} \subset \mathcal{F}^Y$ and (4.11) holds for all $f \in \mathcal{F}$. Moreover, as in the discussion at (2.18), we have $d\mathbb{P}_x|_{\mathcal{M}_t}/d\widehat{\mathbb{P}}_x|_{\mathcal{M}_t} = \widehat{Z}_t$, where $\widehat{Z}$ is the exponential local martingale MF of $Y$ determined by the purely

discontinuous local martingale $\widehat{M}$ with $\widehat{M}_t - \widehat{M}_{t-} = \widehat{\varphi}(Y_{t-}, Y_t)$, and $\widehat{\varphi} := -\varphi/(1+\varphi)$. In short, $X$ can be recovered from $Y$ by a pure-jump Girsanov transformation of the same type that led from $X$ to $Y$. A second application of Theorem 4.4, in which the roles of $X$ and $Y$ are reversed, shows that $\mathcal{F}^Y \subset \mathcal{F}$. $\square$

We now suppose that

$$(4.12) \qquad \int_0^t N(|\varphi|)(X_s)\,dH_s < \infty \qquad \forall\, t \in [0, \zeta[,\ \mathbb{P}_m\text{-a.s.}$$

Thus

$$M_t = \sum_{s \le t} \varphi(X_{s-}, X_s) - (N\varphi * H)_t,$$

the infinite series converging absolutely for each $t \in [0, \zeta[$, $\mathbb{P}_m$-a.s. Consequently, the local martingale AF $M$ has paths of locally finite variation. In particular, $Z$ can be expressed as

$$(4.13) \qquad Z_t = \prod_{s \le t}(1 + \varphi(X_{s-}, X_s))\exp(-(N\varphi * H)_t) \qquad \forall\, t \in [0, \zeta[,$$

the infinite product being absolutely convergent.

Conversely, if $\phi \colon E \times E \to\, ]-1, \infty[$ is symmetric, $\phi(x, \Delta) := 0$ and (4.12) holds, then $Z_t$ in (4.13) defines a positive local martingale MF. Under the family of measures $\mathbb{Q} = (\mathbb{Q}_x \colon x \in E)$ defined by $d\mathbb{Q}_x/d\mathbb{P}_x|_{\mathcal{F}_t} = Z_t$, the process $Y := (X, \mathbb{Q})$ is an $m$-symmetric Markov process whose law is locally absolutely continuous with respect to that of $X$. So we can just start with such a $\phi$ and construct the symmetric process $Y$ in this way. We will now identify the Dirichlet space of $Y$.

Let $\varphi^+ := \varphi \vee 0$ and $\varphi^- := (-\varphi) \vee 0$. Then $0 \le \varphi^- < 1$, $\varphi = \varphi^+ - \varphi^-$ and $1 + \varphi = (1 + \varphi^+)(1 - \varphi^-)$ on $E \times E$. Define, for $t \in [0, \zeta[$,

$$(4.14) \qquad \begin{aligned} Z_t^+ &:= \prod_{s \le t}(1 + \varphi^+(X_{s-}, X_s))\exp\left(\int_0^t N\varphi^-(X_s)\,dH_s\right), \\ Z_t^- &:= \prod_{s \le t}(1 - \varphi^-(X_{s-}, X_s))\exp\left(-\int_0^t N\varphi^+(X_s)\,dH_s\right). \end{aligned}$$

Clearly $Z^+$ is increasing and $Z^-$ is decreasing. Both $Z^+$ and $Z^-$ are MFs that are finite and strictly positive on $[\![0, \zeta[\![$ and $Z = Z^+ \cdot Z^-$. Let $W = (W_t, \mathbb{P}_x^W)$ be the subprocess of $(X, Z^-)$ ("$X$ killed via the MF $Z^-$"). It is easy to see that $W$ coincides with the subprocess of $(Y, 1/Z^+)$.

The AF $\frac{1}{2}\sum_{s \le t} \varphi(X_{s-}, X_s)$ is of bounded variation on compact subintervals of $[0, \zeta[$. We write $\nu := \varphi \cdot J$ for its bivariate Revuz measure, and define

$\mu := 2\nu(1 \otimes \cdot)$. As before,

(4.15)
$$\mathcal{L}^2(\nu) := \left\{ u \in L^2(m) \colon u \text{ has a quasi-continuous version } \tilde{u} \right.$$
$$\left. \text{such that } \int_{E \times E} \overline{u}(x,y)^2 \, \nu(dx,dy) < \infty \right\}.$$

Of course, each element of $\mathcal{F}$ has a quasi-continuous version. Hence we can define

(4.16)
$$\mathcal{F}^\nu := \mathcal{F} \cap \mathcal{L}^2(\nu) = \left\{ u \in \mathcal{F} \colon \int_{E \times E} \overline{u}(x,y)^2 \, \nu(dx,dy) < \infty \right\},$$
$$\mathcal{E}^\nu(u,u) := \mathcal{E}(u,u) + \int_{E \times E} \overline{u}(x,y)^2 \nu(dx,dy), \qquad u \in \mathcal{F}^\nu.$$

By Fatou's lemma, $(\mathcal{E}^\nu, \mathcal{F}^\nu)$ is a symmetric closable quadratic form on $L^2(m)$ if $\varphi \geq c_0 > -1$.

THEOREM 4.6. (i) *Let* $\varphi \colon E \times E \to ]{-}1, \infty[$ *be a symmetric Borel function such that* (4.12) *holds. If* $\varphi(x,y) \geq 0$ *for all* $x, y \in E$, *then*
$$\mathcal{F} \cap L^2(\mu) \subset \mathcal{F}^Y \subset \mathcal{F}^\nu,$$
*where* $\mathcal{F} \cap L^2(\mu)$ *is dense in* $\mathcal{F}^Y$ *with respect to the* $\mathcal{E}_1^Y$*-norm, and for* $u \in \mathcal{F}^Y$,

(4.17)
$$\mathcal{E}^Y(u,u) = \mathcal{E}(u,u) + \int_{E \times E} (u(y) - u(x))^2 \nu(dx,dy).$$

(ii) *If, in addition,* $\varphi$ *is* $J$-*integrable, then* $\mathcal{F}^Y = \mathcal{F}^\nu$.

PROOF. (i) Using the notation in (4.14), we now have (in view of the non-negativity of $\varphi$)
$$Z_t^+ = \prod_{s \leq t} (1 + \varphi(X_{s-}, X_s)), \qquad Z_t^- = \exp\left( -\int_0^t N\varphi(X_s) \, dH_s \right),$$
and $\mu$ is the Revuz measure of $N\varphi * H$. Hence the Dirichlet form $(\mathcal{E}^W, \mathcal{F}^W)$ associated with $W$ and $m$ is given by

(4.18)
$$\mathcal{F}^W = \mathcal{F} \cap L^2(\mu),$$
$$\mathcal{E}^W(u,u) = \mathcal{E}(u,u) + \mu(u^2), \qquad u \in \mathcal{F}^W.$$

But $W$ is also the subprocess $(Y, 1/Z^+)$ of $Y$, so by (4.8) the bivariate Revuz measure of $1/Z^+$ computed with respect to $Y$ and $m$ is $\varphi(1+\varphi)^{-1} \cdot J^Y$, which is nothing but $\nu$. It follows from [27], Theorem II.3.10, that

(4.19)
$$\mathcal{F}^W = \mathcal{F}^Y \cap L^2(\nu),$$
$$\mathcal{E}^W(u,u) = \mathcal{E}^Y(u,u) + \nu(u \otimes u), \qquad u \in \mathcal{F}^W.$$

Combining (4.18) and (4.19), and noting that $\mu$ is smooth with respect to both $X$ and $Y$, we find that $\mathcal{F} \cap L^2(\mu)$ is contained in $\mathcal{F}^Y$ and is dense with respect to the $\mathcal{E}_1^Y$-norm, and for $u \in \mathcal{F} \cap L^2(\mu)$, (4.17) holds.

Assume that $u \in \mathcal{F}^Y$. We may choose a sequence $\{u_n\} \subset \mathcal{F} \cap L^2(\mu)$ such that $u_n \to u$ in $\mathcal{E}_1^Y$-norm. Then $\{u_n\}$ is an $\mathcal{E}_1^Y$-Cauchy sequence and by the result above it is also an $\mathcal{E}_1$-Cauchy sequence. Therefore $u \in \mathcal{F}$ and $u_n \to u$ in $\mathcal{E}_1$-norm and quasi-everywhere (at least along a suitable subsequence). Invoking Fatou's lemma, we have

$$\mathcal{E}(u,u) + \nu(\bar{u}^2) \leq \lim_n (\mathcal{E}(u_n, u_n) + \nu(\bar{u}_n^2)) = \mathcal{E}^Y(u, u) < \infty.$$

It follows that $\mathcal{F}^Y \subset \mathcal{F}^\nu$. As $\bigcup_{k=1}^\infty \mathcal{F}_{F_k} \cap L^2(\mu)$ is $\mathcal{E} + (\cdot, \cdot)_{L^2(\mu)}$-dense in $\mathcal{F} \cap L^2(\mu)$, by Theorem 4.4, (4.17) holds for $u \in \mathcal{F}^Y$. (Here the nest $\{F_k\}$ is as in the statement of Theorem 4.4.)

(ii) We now assume that $\varphi$ is $J$-integrable. For any $u \in \mathcal{F}^\nu$, set $u_n := (u \wedge n) \wedge (-n)$. Then the $J$-integrability of $\varphi$ guarantees $u_n \in \mathcal{F} \cap L^2(\mu)$ and $u_n \to u$ in $\mathcal{E}_1$-norm and q.e. Since $|\bar{u}_n| \leq |\bar{u}|$, we may appeal the dominated convergence theorem and get $u_n \to u$ in $\mathcal{E}_1^\nu$-norm. This implies that $\{u_n\} \subset \mathcal{F}^Y$ is an $\mathcal{E}_1^Y$-Cauchy sequence. Therefore $u \in \mathcal{F}^Y$.   $\square$

In the more general case where (4.12) holds but $\varphi$ is not necessarily positive, we have a weaker result by a similar approach. Let $\nu^+ := \varphi^+ \cdot J$, $\nu^- := \varphi^- \cdot J$, and let $\mu^+ = \nu^+(1 \otimes \cdot)$, $\mu^- = \nu^-(1 \otimes \cdot)$ be the second marginal measures of $\nu^+$ and $\nu^-$, respectively. By [27], Theorem I.4.6, the bivariate Revuz measure of $Z^-$, computed with respect to $X$ and $m$, is

$$\nu_{Z^-}(dx, dy) = \varphi^-(x, y) J(dx, dy) + N\varphi^+(x)\mu_H(dx)\delta_{\{x\}}(dy)$$
$$= \nu^-(dx, dy) + \mu^+(dx)\delta_{\{x\}}(dy)$$

and the bivariate Revuz measure of $1/Z^+$ computed with respect to $(Y, m)$ is

$$\nu_{1/Z^+}^Y(dx, dy) = \frac{\varphi^+}{1 + \varphi^+} J^Y(dx, dy) + N\varphi^-(x)\mu_H(dx)\delta_{\{x\}}(dy)$$
$$= \varphi^+(1 - \varphi^-)J(dx, dy) + \mu^-(dx)\delta_{\{x\}}(dy)$$
$$= \nu^+(dx, dy) + \mu^-(dx)\delta_{\{x\}}(dy).$$

It is now clear that

$$\nu_{Z^-}(1 \otimes \cdot) = \nu_{1/Z^+}^Y(1 \otimes \cdot) = |\mu|.$$

Hence we have by [27], Theorem II.3.10,

$$(4.20) \qquad \begin{aligned} \mathcal{F}^W &= \mathcal{F} \cap L^2(|\mu|), \\ \mathcal{E}^W(u, u) &= \mathcal{E}(u, u) + \nu_{Z^-}(u \otimes u), \qquad u \in \mathcal{F}^W, \end{aligned}$$

and also

$$(4.21) \qquad \begin{aligned} \mathcal{F}^W &= \mathcal{F}^Y \cap L^2(|\mu|), \\ \mathcal{E}^W(u,u) &= \mathcal{E}^Y(u,u) + \nu_{1/Z^+}^Y(u \otimes u), \qquad u \in \mathcal{F}^W. \end{aligned}$$

Combining (4.20) and (4.21), we have the following theorem.

THEOREM 4.7. *Let $\varphi : E \times E \to\ ]-1, \infty[$ be a symmetric Borel function such that (4.12) holds. Then $\mathcal{F} \cap L^2(|\mu|)$ is densely contained in $\mathcal{F}^Y$ and, for $u \in \mathcal{F} \cap L^2(|\mu|)$, we have*

$$(4.22) \qquad \mathcal{E}^Y(u,u) = \mathcal{E}(u,u) + \tfrac{1}{2} \int_0^t (u(y) - u(x))^2 \nu(dx, dy).$$

*Moreover, if $\varphi^+ = 0$, then $\mathcal{F} \subset \mathcal{F}^Y \cap \mathcal{L}^2(\varphi^- \cdot J)$ and (4.22) holds for all $u \in \mathcal{F}$.*

PROOF. It only remains to show the last assertion. In this case, as $\varphi^+ = 0$,

$$(4.23) \quad \frac{1}{Z_t} = \prod_{s \le t} \left(1 + \frac{\varphi^-(Y_{s-}, Y_s)}{1 - \varphi^-(Y_{s-}, Y_s)}\right) \exp\left(-\int_0^t N^Y \frac{\varphi^-}{1 - \varphi^-}(Y_s)\,dH_s\right).$$

But $\varphi^-/(1-\varphi^-) = \varphi^-/(1+\varphi) \ge 0$, so the last assertion follows from Theorem 4.5 with the roles of $X$ and $Y$ interchanged. □

**Acknowledgment.** We gratefully acknowledge the comments of an anonymous referee, whose careful reading of the manuscript helped to improve the exposition of this paper.

Z.-Q. CHEN
DEPARTMENT OF MATHEMATICS
UNIVERSITY OF WASHINGTON
SEATTLE, WASHINGTON 98195
USA
E-MAIL: zchen@math.washington.edu

J. FITZSIMMONS
DEPARTMENT OF MATHEMATICS
UNIVERSITY OF CALIFORNIA, SAN DIEGO
LA JOLLA, CALIFORNIA 92093-0112
USA
E-MAIL: pfitz@euclid.ucsd.edu

M. TAKEDA
MATHEMATICAL INSTITUTE
TOHOKU UNIVERSITY
SENDAI 980-8578
JAPAN
E-MAIL: takeda@math.tohoku.ac.jp

J. YING
INSTITUTE OF MATHEMATICS
FUDAN UNIVERSITY
 SHANGHAI 200433
CHINA
E-MAIL: jying@fudan.edu.cn

T.-S. ZHANG
DEPARTMENT OF MATHEMATICS
UNIVERSITY OF MANCHESTER
OXFORD ROAD, MANCHESTER M13 9PL
UNITED KINGDOM
E-MAIL: tzhang@maths.man.ac.uk